\documentclass[a4paper,10pt]{article}

\usepackage{amsthm}
\usepackage{amsmath}
\usepackage{amssymb}

\usepackage{tikz}
\usepackage{hyperref}
\usepackage{setspace}
\usepackage{mathtools}

\newcommand{\Con}{\ensuremath{\mathcal{C}}}

\newcommand{\D}{\ensuremath{{\mathcal D}}}

\newcommand{\loc}{\ensuremath{\text{loc}}}

\newcommand{\mb}[1]{\ensuremath{\mathbb{#1}}}
\newcommand{\N}{\mb{N}}

\newcommand{\R}{\mb{R}}

\newfont{\bl}{msbm10 scaled \magstep2}

\newcommand{\beq}{\begin{equation}}
\newcommand{\eeq}{\end{equation}}

\newcommand{\notmid}{\mid\kern-0.5em\not\kern0.5em}

\newcommand{\eps}{\varepsilon}

\newenvironment{pr}{\begin{proof}[\textbf{Proof:}] \ }{\end{proof}}
\newtheorem{thm}{Theorem}[section]
\newtheorem{lem}[thm]{Lemma}
\newtheorem{prop}[thm]{Proposition}

\newtheorem{cor}[thm]{Corollary}
\newtheorem{rem}[thm]{Remark}
\newtheorem{defi}[thm]{Definition}

\newcommand{\ep}{\epsilon}

\newcommand{\ttc}{\tilde{\tau}_{\mathrm{co}}}

\newcommand{\Cpq}{C(p,q)}
\newcommand{\Cpqt}{\tilde{C}(p,q)}
\newcommand{\Ch}{C_h(p,q)}
\newcommand{\Lip}{\mathrm{Lip}}

\numberwithin{equation}{section}

\title{Global hyperbolicity for spacetimes with continuous metrics}
\author{Clemens S\"amann\thanks{{\tt clemens.saemann@univie.ac.at}, Faculty of Mathematics, University of Vienna,
Austria.}
}
\date{{November 19, 2019}}

\begin{document}

\maketitle

\begin{abstract}
We show that the definition of global hyperbolicity in terms of the compactness of the causal diamonds and
non-total imprisonment can be extended to spacetimes with continuous metrics, while retaining all of the equivalences to
other notions of global hyperbolicity. In fact, global hyperbolicity is equivalent to the compactness of the space of causal
curves and to the existence of a Cauchy hypersurface. Furthermore, global hyperbolicity implies causal simplicity, stable
causality and the existence of maximal curves connecting any two causally related points.

\vskip 1em

\medskip

\end{abstract}


\section{Introduction}
Global hyperbolicity is the strongest commonly used causality condition in general relativity. It ensures well-posedness of
the Cauchy problem for the wave equation (\cite[Theorem 3.2.11, p.\ 84ff]{BGP:07}), globally hyperbolic spacetimes are the
class of spacetimes used in the initial value formulation of Einstein's equations (\cite[Theorem 16.6, p.\ 177ff.]{Rin:09})
and it plays an important role in the singularity theorems (\cite[Theorem 14.55A, p.\ 431ff. and Theorem 14.61, p.\
436f.]{ONe:83}). These examples emphasize the importance of this notion in Lorentzian geometry and, in particular, in
the theory of general relativity.

Classically (i.e., with $\Con^2$-metric), there are four equivalent notions of global hyperbolicity. These are 
(cf.\ \cite[Subsection 3.11, p.\ 340ff.]{MS:08})
\begin{enumerate}
 \item compactness of the causal diamonds and (strong) causality,
 \item compactness of the space of causal curves connecting two points and causality,
 \item existence of a Cauchy hypersurface,
 \item the metric splitting of the spacetime (cf.\ \cite{BS:03, BS:05}).
\end{enumerate}

General relativity as a geometric theory has been developed for metrics which are smooth (but for all practical purposes 
$\Con^2$ is enough, see e.g.\ \cite{Chr:11, MS:08}), but the PDE point-of-view demands lower regularity in general. In fact, 
even the standard local existence result for the vacuum Einstein equations (\cite{Ren:05}) assumes the metric to be of 
Sobolev-regularity $H^s_\loc$ (with $s>\frac{5}{2}$). Furthermore, in recent years the regularity of the metric has been 
lowered even more (e.g.\ \cite{KRS:12}). The critical regularity class, where many aspects of causality theory work as in 
the smooth case is $\Con^{1,1}$ (\cite{CG:12, Min:13, KSS:14, KSSV:14}). Below $\Con^{1,1}$ the geodesic 
equations need not be uniquely solvable (\cite{HW:51}) and so the exponential map cannot be used to locally transport 
notions from the tangent space to the manifold, which is an indispensable tool in semi-Riemannian geometry.

Not only the PDE point-of-view demands low regularity, also, physically relevant models of spacetimes impose certain 
restrictions on the regularity of the metric. In particular, modeling different types of matter in a spacetime might lead to 
a discontinuous energy-momentum tensor, and hence by the Einstein equations to metrics of regularity below 
$\Con^2$ (\cite{Lic:55, MS:93}). Prominent examples are spacetimes which model the inside and outside of a star or shock 
waves. There are even less regular, physically interesting models like spacetimes with conical singularities and cosmic 
strings (\cite{Vic:90, VW:00}), (impulsive) gravitational waves (\cite{Pen:72a, KS:99, SS:12}, especially \cite{LSS:14, 
PSSS:14}, where the Lipschitz continuous form of the metric is used) and ultrarelativistic black holes (e.g.\ \cite{AS:71}).

Motivated by providing some tools for studying the Cauchy problem in low regularity (\cite{Chr:13}), Chru\'sciel and Grant 
developed causality theory for spacetimes with continuous metrics (\cite{CG:12}). They showed how to approximate continuous 
metrics by smooth ones while retaining a control of the causality properties of the approximating metrics. Using these 
methods they proved that domains of dependence are globally hyperbolic and admit smooth time functions, for example. 
Moreover, some properties cease to hold (as one would expect from the classical Riemannian examples by Hartman and Wintner, 
\cite{HW:51}) E.g.\ the push-up principle is not valid anymore and lightcones could ``bubble up'' (i.e., they are not 
hypersurfaces).
\bigskip

There have been several approaches to global hyperbolicity for non-smooth metrics. The approach of Clarke (\cite{Cla:98})
advocated the view that singularities should not be understood as the obstruction to extend geodesics but as obstructions to
the evolution of physical fields. Vickers and Wilson (\cite{VW:00}) showed that spacetimes with conical
singularities and cosmic strings are globally hyperbolic in this sense.

Sorkin and Woolgar (\cite{SW:96}) used order-theoretic methods to define and investigate the notion of $K$-causality, which 
is the smallest relation containing the usual timelike relation $I^\pm$ that is transitive and closed. Based on this 
relation, they define a notion of global hyperbolicity in terms of the compactness of the $K$-causal diamonds. Moreover, 
their concept of global hyperbolicity agrees with the usual notion for metrics that are $\Con^2$.

Fathi and Siconolfi (\cite{FS:12}) investigated the existence of smooth time functions in the setting of \emph{cone 
structures}. These are families of certain closed, convex cones in the tangent spaces of a manifold, which may but need 
not arise from the (forward) lightcones of a Lorentzian metric. Then they establish a notion of global 
hyperbolicity which assumes compactness of the causal diamonds and stable causality (this is stronger than the usual 
definition for spacetimes). Their major result is that for continuous and globally hyperbolic cone structures there is a 
smooth (Cauchy) time function. Moreover they also show that global hyperbolicity is stable, i.e., for every continuous and 
globally hyperbolic cone structure there is a globally hyperbolic cone structure with wider cones.

H\"ormann, Kunzinger and Steinbauer (\cite[Definition 6.1, p.\ 182]{HKS:12}) defined global hyperbolicity in terms of the
metric splitting, which is well-suited for the Cauchy problem for the wave equation on non-smooth spacetimes (with
weakly-singular metrics). Furthermore in \cite{HS:14} this notion was investigated for a class of non-smooth wave-type
spacetimes (generalizations of pp-waves, with non-flat wave surfaces).
\bigskip 

Building on methods developed in \cite{CG:12} we propose a notion of global hyperbolicity that extends the usual 
notion of global hyperbolicity based on the compactness of the causal diamonds and still retains all equivalences as above 
(the metric splitting has to be weakened to a topological splitting of course). This, in particular, shows that global 
hyperbolicity is still the strongest of the commonly used causality conditions.
\bigskip

The outline of the article is as follows: In the remainder of the introduction we fix notation and state results used 
throughout the article. Then in Section \ref{sec-top} we define topologies on spaces of causal curves connecting 
two points and investigate their relationship. Results and methods of this section will then be used for the definition of 
global hyperbolicity and the equivalence to the compactness of the space of causal curves in Section \ref{sec-gh}. 
Furthermore, in Section \ref{sec-stab} we will show that global hyperbolicity implies stable causality and give a 
self-contained proof of the stability of global hyperbolicity. This will allow us to show the equivalence of global 
hyperbolicity to the existence of a Cauchy hypersurface in Section  \ref{sec-chs}. Finally, in Section \ref{sec-max}, we 
apply our results to obtain existence of maximal curves connecting any two causally related points in a globally hyperbolic 
spacetime.

\subsection*{Preliminaries}
Let $M$ be a smooth connected manifold and $g$ a continuous Lorentzian metric. Moreover let $(M,g)$ be
time-oriented (i.e., there is a continuous timelike vector field), fix a smooth complete Riemannian metric $h$ on $M$ and
denote the induced metric by $d^h$. 

\begin{defi}
 Let $\lambda : I \rightarrow M$ be a locally Lipschitz continuous curve (with respect to $d^h$; in fact, it is independent
of the choice of complete Riemannian metric, see \cite[Proposition 2.3.1, p.\ 14]{Chr:11}), then it is called (cf.\
\cite[Definition 1.3 p.\ 4]{CG:12})
\begin{enumerate}
 \item \emph{timelike} if $g(\dot\lambda,\dot\lambda)<0$ almost everywhere,
 \item \emph{causal} if $g(\dot\lambda,\dot\lambda)\leq 0$ and $\dot\lambda\neq0$ almost everywhere.
\end{enumerate}
A causal curve $\lambda$ is called \emph{future (past) directed} if $\dot\lambda$ belongs to the future (past) lightcone
almost everywhere.
\end{defi}

For $p,q\in M$ we define the space of future directed causal curves connecting $p$ and $q$ by
\begin{equation*}
\Cpq:=\{\lambda:I_\lambda\rightarrow M : \lambda \text{ future directed causal curve from }p\text{ to }q\}\ \slash \sim \,,
\end{equation*}
where $I_\lambda\subseteq \R$ denotes a compact interval. Moreover, $\lambda\sim\gamma$ if there exists an orientation 
preserving reparametrization, i.e., a map $\phi\colon I_\lambda\rightarrow I_\gamma$ with $\lambda = \gamma \circ 
\phi$ that is absolutely continuous, surjective, strictly monotonically increasing and its inverse is absolutely continuous. 
So we consider two curves equal if one is a reparametrization of the other (while keeping the orientation). Note that in this 
regularity class it is not clear that one can make this identification, thus we clarify this in the following lemma.

\begin{lem}\label{lem-cpq-equiv}
 \leavevmode
 \begin{enumerate}
  \item The relation $\sim$ is an equivalence relation.
  \item \label{lem-uniq-param} For every class $\lambda\in\Cpq$ there is a unique parametrization with respect to 
$h$-arclength. Moreover, there is a unique parametrization on $[0,1]$ proportional to $h$-arclength.
 \end{enumerate}
\end{lem}
\begin{pr}
 \begin{enumerate}
  \item It remains to show that $\sim$ is transitive because the composition of two absolutely continuous functions need not 
be absolutely continuous in general. However, in our case the composition of two parametrizations as above is absolutely 
continuous. Let $f\colon[a,b]\rightarrow [r,s]$, $g\colon[r,s]\rightarrow[c,d]$ be two parametrizations as above and set 
$h:=g \circ f\colon [a,b]\rightarrow [c,d]$. Then $h$ is a strictly monotonically increasing homeomorphism and thus has 
bounded variation. This implies that $h$ is absolutely continuous by \cite[Theorem IX.3.5, p.\ 252]{Nat:55}.
 
 \item Let $\lambda\in\Cpq$ and $\tilde{\lambda}$ a representative of $\lambda$ defined on $[a,b]$. It is well-known (cf.\ 
\cite[p.\ 3]{Min:08a}) that $\tilde{\lambda}$ has a parametrization with respect to $h$-arclength, i.e., 
$h(\dot\lambda,\dot\lambda)=1$ almost everywhere. This parametrization is given by $\phi\colon[a,b]\rightarrow[0,L]$,
$\phi(t):=L^h(\tilde{\lambda}|_{[0,t]})$ $(t\in[a,b])$, where $L^h$ denotes the length functional with respect to the 
Riemannian metric $h$ and $L:=L^h(\lambda)$. Clearly, $\phi$ is strictly monotonically increasing, surjective and Lipschitz 
continuous (with Lipschitz constant bounded by the Lipschitz constant of $\tilde{\lambda}$). Thus $\phi$ is absolutely 
continuous. It remains to show that $\phi^{-1}$ is absolutely continuous. Note that $\dot\phi>0$ almost everywhere, so 
$\{t\in[a,b]: \dot\phi(t)=0\}$ has measure zero. Then, by a result of Zareckii (\cite[p.\ 271]{Nat:55}), $\phi^{-1}$ is 
absolutely continuous. Consequently $\tilde{\lambda}\sim \tilde{\lambda}\circ \phi^{-1}$, and $\tilde{\lambda}\circ\phi^{-1}$ 
is parametrized with respect to $h$-arclength.

Set $\lambda^*:=\tilde{\lambda}\circ\phi^{-1}$ and assume that we have another parametrization with respect to 
$h$-arclength, i.e.,  $\lambda^*\circ\psi\colon[0,L]\rightarrow M$. Then almost everywhere $1=h((\lambda^*\circ\psi)\dot\ 
,(\lambda^*\circ\psi)\dot\ ) = (\dot{\psi})^2 h(\dot{\lambda^*},\dot{\lambda^*})= (\dot\psi)^2$ and since $\dot\psi\geq 0$ 
we conclude that $\dot\psi = 1$. So $\psi$ is the identity on $[0,L]$ (because $\psi$ is absolutely continuous).

Finally, reparametrizing $\lambda^*$ by $s\mapsto L^h(\lambda)s$ gives a parametrization on $[0,1]$ 
proportional to $h$-arclength, which is also unique. Note that then the best Lipschitz constant (with respect to $d^h$) of 
$\lambda$, denoted by $\Lip^h(\lambda)$, is equal to $L^h(\lambda)$.
 \end{enumerate}
\end{pr}

Now that we fixed this terminology, we will from now on synonymously say curve for a fixed parametrization and its class, if 
it is clear from the context which one is meant. Moreover, we will denote the images of curves always by the corresponding 
capital letters, e.g.\ $\Lambda:=\lambda(I_\lambda)$.
\bigskip

We need that we can approximate a continuous metric arbitrarily well by smooth metrics and so we introduce the distance 
between two metrics $g$ and $g'$ by
\begin{equation*}
\Delta(g,g'):=  \sup_{p\in M,\, 0\neq X,Y\in T_p M} \frac{|g(X,Y)-g'(X,Y)|}{|X|_h |Y|_h}\,.
\end{equation*}

\begin{prop}\label{prop-cg-2.1} (\cite[Proposition 2.1, p.\ 3]{CG:12})
For every $\ep>0$ there are smooth Lorentzian metrics $\hat{g}$ and $\check{g}$, such that $\check{g}\prec g\prec\hat{g}$ and
\begin{equation*}
 \Delta(g,\check{g}) + \Delta(g,\hat{g}) \leq \ep\,.
\end{equation*}
Moreover, $g\prec\hat{g}$ means that the lightcones of $\hat{g}$ are strictly greater than those of $g$, i.e., if a
non-zero vector is causal for $g$ then it is timelike for $\hat{g}$. Additionally we mean by $g_1\preceq g_2$ that every 
$g_1$-causal vector is causal for $g_2$. 
\end{prop}
Furthermore, there is a sequence of smooth metrics that converges locally uniformly to $g$ and this sequence can be chosen 
to be monotonically increasing or decreasing (\cite[Proposition 2.3, p.\ 5f.]{KSV:15}). This implies that on compact sets we
can always have the sequence smaller than a given metric. 

\begin{lem}\label{lem-small-metri}
 Let $g$ be a continuous metric, $g\prec g_{n+1}\prec g_n$ for all $n\in\N$, where $(g_n)_n$ converges locally uniformly to 
$g$. Moreover let $K\subseteq M$ be compact and $g'\succ g$, then there is an $n_0\in\N$ such that for all $n\geq n_0$ we 
have $g_n\prec g'$ on $K$.
\end{lem}
\begin{pr}
 Set $L:=\{X\in TM|_K: |X|_h=1,\ g(X,X)\leq 0\}$, then $L$ is a compact subset of $TM$. We claim that there is $\delta>0$ 
such that $\{X\in TM|_K: |X|_h=~1,\  g(X,X)<\delta\}\subseteq \{X\in TM: g'(X,X)<0\}$. Assume to the contrary that there 
are $X_k\in TM|_K$, $|X_k|_h=1$, $g(X_k,X_k)<\frac{1}{k}$ but $g'(X_k,X_k)\geq 0$ for $k\in\N$. By compactness we can 
without 
loss of generality assume that $X_k\to X\in TM|_K$ with $|X|_h=1$ and $g(X,X)\leq 0$ but $g'(X,X)\geq 0$ --- a contradiction 
to $g\prec g'$. Choose $n_0\in\N$ such that for all $n\geq n_0$ we have $\Delta(g_n,g)< \delta$. Let $n\geq n_0$, $X\in 
TM|_K$, $|X|_h=1$ and $g_n(X,X)\leq 0$. Then $g(X,X)\leq g_n(X,X) + \Delta(g_n,g)<\delta$, hence $g'(X,X)<0$. Thus $g_n\prec 
g'$ on $K$.
\end{pr}

We need the following limit curve theorem from $\cite{CG:12}$ for continuous metrics (see also \cite[Theorem 3.1,(1), p.\ 
8f.]{Min:08a}), which we slightly strengthen.
\begin{thm}\label{thm-cq-1.6}(\cite[Theorem 1.6, p.\ 6]{CG:12})
Let $(g_n)_n$ be a sequence of smooth metrics such that $g\prec g_{n+1}\prec g_n$ for all $n\in\N$ and $g_n\to g$ locally 
uniformly. Let $(\gamma_n)_n$ be a sequence of (parametrized) curves accumulating at some $p\in M$ such that $\gamma_n$ is 
$g_n$-causal. Then, if
\begin{enumerate}
 \item the $\gamma_n$'s are all defined on the same interval, say $[a,b]$, and have uniformly bounded Lipschitz constants, or
 \item the $\gamma_n$'s are inextendible,
\end{enumerate}
there exists a causal curve $\gamma$ through $p$ and there is a subsequence $(\gamma_{n_k})_k$ of $(\gamma_n)_n$ which 
converges to $\gamma$ uniformly on compact sets. In particular, this implies in the first case uniform convergence on 
$[a,b]$ and in the second case that $\gamma$ is inextendible.
\end{thm}
\begin{pr}
 As in the proof of \cite[Theorem 1.6, p.\ 6]{CG:12} all $\gamma_n$ are $g_0$-causal, hence by the smooth result 
(\cite[Theorem 3.1,(1), p.\ 8f.]{Min:08a}), there is a $g_0$-causal limit curve and a subsequence $(\gamma_{n_k})_k$ of 
$(\gamma_n)_n$ such that $\gamma_{n_k}\to \gamma$ locally uniformly. In the first case $\gamma$ is defined on $[a,b]$ and it 
the second case it is inextendible. It remains to show that $\gamma$ is $g$-causal. Since this is a local question, we can 
restrict to some finite interval $[c,d]$ ($[a,b]$ in the first case, some finite interval in the second case). By uniform 
convergence on $[c,d]$, all segments $\gamma_n|_{[c,d]}$ are contained in some open, relatively compact open neighborhood 
$U$ of $\gamma([c,d])$. Let $\hat{g}\succ g$, then by Lemma \ref{lem-small-metri} there is a $n_0\in\N$ such that for all 
$n\geq n_0$, $g_n\prec \hat{g}$ on $\overline{U}$. Let $k_0\in\N$ be such that $n_k\geq n_0$ for all $k\geq k_0$, then 
$(\gamma_{n_k}|_{[c,d]})_{k\geq k_0}$ is a sequence of $\hat{g}$-causal curves, which has (as above) a limit curve that is 
$\hat{g}$-causal. However, this curve has to be $\gamma|_{[c,d]}$, since this is a subsequence of a converging sequence.
Summing up, this shows that $\gamma$ is $\hat{g}$-causal for every $\hat{g}\succ g$, hence $\gamma$ is $g$-causal by 
\cite[Proposition 1.5, p.\ 5f.]{CG:12}.
\end{pr}

Concerning chronological and causal futures and pasts, our notation follows mainly \cite{CG:12}. In particular,
$I^\pm_{g'}(A)$, $J^\pm_{g'}(A)$ denotes the set of all points in $M$ which can be reached by a future/past directed
$g'$-timelike respectively $g'$-causal curve from $A\subseteq M$. Moreover $\check{I}^\pm(A):=\bigcup_{\check{g}\prec
g}I^\pm_{\check{g}}(A)$ (where $\check{g}$ is smooth), which is open by \cite[Proposition 1.4, p.\ 5]{CG:12}.

\section{Topologies on spaces of causal curves}\label{sec-top}
In this section we are going to clarify which topology to put on $\Cpq$ and how to handle parametrized curves. In doing so 
we are going to resolve several technicalities and this will simplify the methods later on. We start with the topology 
$\tau$ as introduced in \cite[p.\ 208]{HE:73} (for smooth metrics).

Let $U\subseteq M$ be open, then we set $O(U):= \{\lambda\in \Cpq:\Lambda\subseteq U\}$.
\begin{lem}
The sets $O(U)$, for $U\subseteq M$ open, form a basis for a topology on $\Cpq$, which we call $\tau$.
\end{lem}
\begin{pr}
Obviously $O(U_1 \cap U_2) = O(U_1)\cap O(U_2)$ for $U_1,U_2\subseteq M$ open.
\end{pr}
By construction, a set $\mathcal{U}\subseteq C(p,q)$ is $\tau$-open if it can be written as a union 
of sets $O(U)$, i.e., $\mathcal{U}=\bigcup_{\alpha\in A} O(U_\alpha)$, where $U_\alpha\subseteq M$ is open in $M$ for 
$\alpha\in A$.
\bigskip

There is a connection between topological properties of $\tau$ and causality of $(M,g)$, as the following result shows.
\begin{lem}
Let $(M,g)$ be a continuous spacetime.
\begin{enumerate}
 \item If $(M,g)$ is causal and $\lambda,\gamma\in \Cpq$ with $\Lambda\subseteq\Gamma$, then $\lambda = \gamma$.
 \item $(M,g)$ is causal if and only if the topology $\tau$ is Hausdorff on $\Cpq$  for all $p,q\in M$.
\end{enumerate}
\end{lem}
\begin{pr}
\begin{enumerate}
 \item
Let $\lambda,\gamma\in C(p,q)$ with $\Lambda\subseteq \Gamma$. Moreover we can assume that $\lambda,\gamma$ are defined on 
$[0,1]$ with $\lambda(0)=\gamma(0)=p$ and $\lambda(1)=\gamma(1)=q$. Now we define the function $f\colon[0,1]\rightarrow 
[0,1]$ by $f(t):=\inf\{s\in[0,1]: \lambda(t)=\gamma(s)\}$. Then it is clear that $\lambda(t)=\gamma(f(t))$ for all 
$t\in[0,1]$. First we show that $f$ is strictly monotonically increasing: let $t_1<t_2$ be in $[0,1]$ and suppose that 
$f(t_1)\geq f(t_2)$. Then $\lambda(t_1) < \lambda(t_2) = \gamma(f(t_2))\leq \gamma(f(t_1)) = \lambda(t_1)$, where the order 
relations are meant to be causal relations. This gives a contradiction since $(M,g)$ is causal.

At this point we show that $f$ is continuous: suppose $\exists t_0\in[0,1]$ such that $\lim_{t\nearrow
t_0}f(t)<\lim_{t\searrow t_0}f(t)$. It follows that $\lambda(t_0)=\lim_{t\nearrow t_0}\lambda(t) = \lim_{t\nearrow
t_0}\gamma(f(t))< \lim_{t\searrow t_0}\gamma(f(t)) = \lim_{t\searrow t_0}\lambda(t) = \lambda(t_0)$ --- again a contradiction
to the causality of $(M,g)$.

Furthermore, $f$ is surjective. In fact, let $s\in[0,1]$, then $f(0)=0 \leq s \leq 1=f(1)$ and hence by the intermediate 
value theorem there is a $t\in[0,1]$ with $f(t)=s$. This allows us to conclude that $\Lambda=\Gamma$. Assume the 
contrary, i.e., there is a $x\in\Gamma\backslash\Lambda$. Let $s\in[0,1]$ with $\gamma(s)=x$, then since $f$ is surjective, 
there is a $t\in[0,1]$ with $f(t)=s$. Consequently, $x=\gamma(s)=\gamma(f(t))=\lambda(t)$, so $x\in\Lambda$ --- a 
contradiction.

Then by Lemma \ref{lem-cpq-equiv},\ref{lem-uniq-param} there are reparametrizations 
$\phi,\psi:[0,L^h(\lambda)]\rightarrow[0,1]$ such that $\lambda\circ\phi$ and $\gamma\circ\psi$ are parametrized with 
respect to $h$-arclength. This yields $\lambda\sim\gamma$ via $\psi\circ\phi^{-1}$, because 
$\lambda\circ\phi=\gamma\circ\psi$. To see this last claim let $s\in[0,1]$, then there is a $s'\in[0,1]$ with 
$\lambda(\phi(s))=\gamma(\psi(s'))$. As above, the images of the curves $\lambda\circ\phi|_{[0,s]}$ and 
$\gamma\circ\psi|_{[0,s']}$ agree. Then $s = L^h(\lambda\circ\phi|_{[0,s]}) = L^h(\gamma\circ\psi|_{[0,s']}) = s'$, since 
both curves are parametrized with respect to $h$-arclength and thus $\lambda(\phi(s))=\gamma(\psi(s))$. 

\item
Assume $(M,g)$ to be causal and let $\lambda,\gamma\in C(p,q)$ with $\lambda\neq\gamma$, then by the first part
$\Gamma\nsubseteq\Lambda$ and $\Lambda\nsubseteq\Gamma$, so there is an $x\in\Gamma\backslash\Lambda$ and a
$y\in\Lambda\backslash\Gamma$. Since $M$ is Hausdorff there exist open neighborhoods $U_x,V_y$ of $x$ and $y$, respectively, 
such that $U_x\cap V_y=\emptyset$. Assume $\gamma\colon[a,b]\rightarrow M$ and let $t\in[a,b]$ such that $\gamma(t)=x$. Then
setting $\tilde{\Gamma}_1:=\gamma([a,t])$ and $\tilde{\Gamma}_2:=\gamma([t,b])$, the sets
$\Gamma_i:=\tilde{\Gamma}_i\backslash U_x$ are compact too ($i=1,2$). Hence we can separate them by open sets $U_1, U_2$
and setting $U:=U_1\cup U_x\cup U_2$ yields an open neighborhood of $\Gamma$ such that $U\backslash U_x$ is not connected.
Performing an analogous procedure for $\lambda, y$ and $V_y$ one obtains an open neighborhood $V=V_1\cup V_y\cup V_2$ of
$\Lambda$ such that $V\backslash V_y = V_1\cup V_2$ is not connected. Then $O(U)$ (respectively $O(V)$) is an open
neighborhood of $\gamma$ (respectively $\lambda$) such that $O(U)\cap O(V) = \emptyset$. To verify the last claim let 
$\xi\in O(U)\cap O(V)=O(U\cap V)$ and note that $\xi$ has to start in $U_1\cap V_1$ and end in $U_2\cap V_2$. However these 
open sets are disjoint and so $\xi$ has to pass through $U_x\cap V_y=\emptyset$ --- a contradiction.

Assume now that $(M,g)$ is not causal. Then there is a closed causal curve $\gamma$. Choose $p,q\in\Gamma$ (and without loss
of generality $p< q$), then the segment of $\gamma$ connecting $p$ to $q$, denoted by $\lambda$, cannot be separated
in $\Cpq$ from the causal curve $\tilde{\gamma}$, where $\tilde{\gamma}$ is obtained by the loop $\gamma$ at $p$ joined with
$\lambda$. Therefore $(\Cpq,\tau)$ is not Hausdorff.
\end{enumerate}
\end{pr}
Note that also a continuous spacetime cannot be compact and chronological (hence causal). This follows from the
smooth case (cf.\ \cite[Lemma 14.10, p.\ 407]{ONe:83}) by choosing a smooth metric $\check{g}\prec g$, since then there is a
closed $\check{g}$-timelike curve, which is also $g$-timelike. 
\bigskip

From this point on we assume $(M,g)$ to be causal and we define $d$ to be the Hausdorff distance restricted to the images of
curves in $\Cpq$. To be more precise:
\begin{equation*}
d(\lambda,\gamma):=\inf\{r>0: \Lambda \subseteq (\Gamma)_r\text{ and } \Gamma \subseteq (\Lambda)_r \}\,,
\end{equation*}
where by $(A)_r$ we denote $\{x\in M: d^h(A,x)<r\}$ (for $A\subseteq M$). It is indeed a metric on $C(p,q)$ since it is a 
metric on all closed subsets of $M$, $(M,g)$ is causal and we restrict to future directed curves. We denote the metric 
topology of $d$ by $\tau'$.

\begin{lem}\label{lem-tau-leq-tau'}
 The topology $\tau'$ is finer than $\tau$, i.e., $\tau\subseteq\tau'$.
\end{lem}
\begin{pr}
Let $\lambda_0\in C(p,q)$ and let $U\subseteq M$ be an open neighborhood of $\Lambda_0$. Then it suffices to prove
that there is a $\delta>0$ such that $B_\delta^d(\lambda_0) \subseteq O(U)$, where $B_\delta^d(\lambda_0)$ denotes the ball
with respect to the Hausdorff metric $d$ of radius $\delta$ around $\lambda_0$. Since $U$ is open and $\Lambda_0$ is compact,
there is a $\delta>0$ such that $(\Lambda_0)_\delta\subseteq U$. Then let $\lambda\in B_\delta^d(\lambda_0)$, i.e.,
$d(\lambda, \lambda_0)<\delta$, therefore, in particular, $\Lambda\subseteq (\Lambda_0)_\delta$, which is equivalent to
$\lambda\in O((\Lambda_0)_\delta)$.
\end{pr}

By choosing a special parametrization of a curve in $\Cpq$ we can define (cf.\ \cite{Cho:68})
\begin{align*}
 \Ch:=&\{\lambda\colon[0,1]\rightarrow M : \lambda \text{ Lipschitz continuous, future directed, }\\
 &\text{causal, } \lambda(0)=p, \lambda(1)=q \text{ and } h(\dot\lambda,\dot\lambda)=\text{constant a.e.\ }\}\,.
\end{align*}

Note that $\Ch$ is a subset of $\Con([0,1],M)$, which comes with the natural topology of uniform convergence or 
equivalently, the compact-open topology, which we denote by $\ttc$. This allows us to define $\Cpqt$ as the closure of $\Ch$ 
with respect to $\ttc$, i.e.,
\begin{equation*}
 \Cpqt:= \overline{\Ch}^{\ttc}\,.
\end{equation*}

By slight abuse of notation we denote the compact-open topology of $\Con([0,1],M)$ restricted to $\tilde{C}(p,q)$ also by 
$\ttc$. Moreover, the topology $\ttc$ is induced by the metric $\rho$ given by
\begin{equation*}
 \rho(\lambda,\gamma)=\sup_{t\in[0,1]} d^h(\lambda(t),\gamma(t))\quad (\lambda,\gamma\in\Cpqt)\,.
\end{equation*}

\begin{lem}
 The elements of $\Cpqt$ are future directed causal curves from $p$ to $q$ defined on $[0,1]$.
\end{lem}
\begin{pr}
 This follows from Theorem \ref{thm-cq-1.6}.
\end{pr}
Note that a limit of curves in $\Ch$ need not be parametrized proportional to $h$-arclength, hence the need to take the 
closure of $\Ch$ in $\Con([0,1],M)$.
\bigskip

We define the map $\Phi\colon \Cpqt \rightarrow \Cpq$ by assigning to a curve $\tilde{\lambda}\in\Cpqt$ its class 
in $\Cpq$.
\begin{lem}\label{lem-phi-cont}
 The map $\Phi\colon \Cpqt \rightarrow \Cpq$, $\Cpqt\ni\tilde{\lambda}\mapsto [\tilde{\lambda}]$ is surjective and 
$\ttc$-$\tau$ continuous.
\end{lem}
\begin{pr}
 We can always parametrize a curve in $\Cpq$ on $[0,1]$ proportional to $h$-arclength by Lemma 
\ref{lem-cpq-equiv},\ref{lem-uniq-param}, so $\Phi$ is surjective.

To see $\ttc$-$\tau$ continuity, let $\tilde{\gamma}\in\Cpqt$ and set $\gamma:=\Phi(\tilde{\gamma})$. Let $U\subseteq M$ be 
open such that $\gamma\in O(U)$. Then, since $\Gamma$ is a compact subset of $U$, there is a $\delta>0$ such that 
$(\Gamma)_\delta\subseteq U$. We claim that $B^\rho_\delta(\tilde{\gamma})\subseteq \Phi^{-1}(O(U))$. Let 
$\lambda\in B^\rho_\delta(\tilde{\gamma})$ and $t\in[0,1]$, then $d^h(\lambda(t),\Gamma)\leq 
d^h(\lambda(t),\tilde{\gamma}(t))<\delta$. Consequently, $\Lambda\subseteq (\tilde{\Gamma})_\delta = 
(\Gamma)_\delta\subseteq U$. Thus $\Phi^{-1}(O(U))$ is a $\ttc$-open neighborhood of $\tilde{\gamma}$.
\end{pr}
Note that, by Lemma \ref{lem-cpq-equiv},\ref{lem-uniq-param}, $\Phi|_{\Ch}\colon\Ch\rightarrow\Cpq$ is bijective.
A consequence of the above lemma is that $\Cpq$ is $\tau$-compact if $\Cpqt$ is $\ttc$-compact. We also want the converse to 
hold, thus we will show that under a stronger causality condition $\Phi$ is a proper map, i.e., preimages of compact 
sets are compact.
\bigskip

Here we introduce two imprisoning conditions (cf.\ \cite{Min:08b}).
\begin{defi}\label{def-imp}
\leavevmode
 \begin{enumerate}
  \item A spacetime is called \emph{non-totally imprisoning} if there is no future or past inextendible causal curve
contained in a compact set.
  \item A spacetime is called \emph{non-partially imprisoning} if there is no future or past inextendible causal curve
returning infinitely often to a compact set. 
 \end{enumerate}
\end{defi}
Clearly, non-partially imprisoning implies non-totally imprisoning, which implies causality. Moreover, strong causality 
implies non-partially imprisoning, since the proof of \cite[Lemma 14.13, p.\ 408]{ONe:83} works also for continuous metrics.

The next proof was suggested to us by E.\ Minguzzi.
\begin{lem}\label{lem-unif-lip}
 The spacetime $(M,g)$ is non-totally imprisoning if and only if for every $K\subseteq M$ compact there exists a 
$C>0$ such that the $h$-arclengths of all causal curves contained in $K$ are bounded by $C$.
\end{lem}
\begin{pr}
Since an inextendible causal curve has infinite $h$-arclength, the converse is clear.

Assume that $(M,g)$ is non-totally imprisoning but there is no bound on the $h$-arclengths of causal curves contained in 
some compact set $K$, i.e., there is a sequence $(\gamma_n)_n$ of causal curves contained in $K$ such that  
$L^h(\gamma_n)\to\infty$ as $n\to\infty$. Without loss of generality we can assume that the sequence is parametrized with 
respect to $h$-arclength, $\gamma_n\colon[0,L^h(\gamma_n)]\rightarrow K$ and since $K$ is compact $\gamma_n(0)\to p\in K$. 
Then as in the proof of Theorem \ref{thm-cq-1.6} the limit curve theorem (\cite[Theorem 3.1, p.\ 8ff.]{Min:08a}) yields the 
existence of an inextendible causal curve contained in $K$ --- a contradiction to $(M,g)$ being non-totally imprisoning.
\end{pr}
A consequence of this lemma is that one obtains a bound on the Lipschitz constants of curves in $\Ch$, if they are contained 
in some compact set and the spacetime is non-totally imprisoning (since for every $\gamma\in\Ch$ one has 
$L^h(\gamma)=\Lip(\gamma)$ ($=\sqrt{h(\dot\gamma,\dot\gamma)}$ almost everywhere)).
\bigskip

At this point we are able to show that $\Phi$ is a proper map, if $(M,g)$ is non-totally imprisoning.
\begin{thm}
Let $(M,g)$ be non-totally imprisoning. Then $\Phi\colon \Cpqt \rightarrow \Cpq$ is $\ttc$-$\tau$ proper. 
\end{thm}
\begin{pr}
Let $K\subseteq \Cpq$ be compact and set $\tilde{K}:=\Phi^{-1}(K)$. Since $\ttc$ is a metric topology, compactness is 
equivalent to sequential compactness (\cite[Theorem 4.1.17, p.\ 256]{Eng:89}) and thus it suffices to consider a sequence
$(\tilde{\gamma}_n)_n$ in $\tilde{K}$. Setting $\gamma_n:=\Phi(\tilde{\gamma}_n)$ yields a sequence in $K$. Thus, by
assumption and since $\tau$ is first countable (\cite[Theorem 3.10.31, p.\ 209]{Eng:89}), there is a subsequence
$(\gamma_{n_k})_k$ of $(\gamma_n)_n$ that converges to some $\gamma\in K$ with respect to $\tau$. Let 
$U\subseteq M$ be an open, relatively compact neighborhood of $\Gamma$, then there is a $k_0\in\N$ such that for all $k\geq 
k_0$: $\Gamma_{n_k}\subseteq U$. This implies, in particular, that the $\tilde{\gamma}_{n_k}$'s are contained in $\bar{U}$. 
Hence by Lemma \ref{lem-unif-lip} there is a bound on the Lipschitz constants of $(\tilde{\gamma}_{n_k})_k$ and since they 
all start at $p$, we can apply Theorem \ref{thm-cq-1.6} to get a $\ttc$-convergent subsequence. This shows that $\tilde{K}$ 
is $\ttc$-compact.
\end{pr}

As mentioned above, we now have the following characterization of compactness of the space of causal curves connecting two 
points.
\begin{cor}\label{cor-phi-proper}
 Let $(M,g)$ be non-totally imprisoning. Then $\Cpq$ is $\tau$-compact if and only if $\Cpqt$ is $\ttc$-compact.
\end{cor}

Because we always want to work with sequences of curves instead of nets, we establish below that $\tau=\tau'$ and hence 
compactness is equivalent to sequential compactness.

\begin{lem}
 Let $(M,g)$ be non-totally imprisoning, then $\tau=\tau'$ and hence $\tau$ is metrizable.
\end{lem}
\begin{pr}
  Assume that $\tau\subsetneq\tau'$, i.e., there is a $\lambda\in\Cpq$ and $\ep>0$ such that for every open
neighborhood $U\subseteq M$ of $\Lambda$ we have $O(U)\not\subseteq B^d_\ep(\lambda)$. For $n\in\N, n\geq 1$ set
$U_n:=(\Lambda)_{\frac{1}{n}}$, then our assumption yields a sequence $(\gamma_n)_n$  with $\gamma_n\in O(U_n)$ and
$d(\gamma_n,\lambda)\geq \ep$ for all $n\in\N, n\geq 1$. Furthermore, it is clear that $\gamma_n\to\lambda$ with respect to
$\tau$. 

For every $n\in\N$, $n\geq 1$ there is a unique $\tilde{\gamma}_n\in \Phi^{-1}(\{\gamma_n\})\cap \Ch$, hence this 
yields a sequence in $\Ch$. This sequence satisfies: $\tilde{\Gamma}_n = \Gamma_n\subseteq U_1$, which is an open, 
relatively compact neighborhood of $\Lambda$. So by Lemma 
\ref{lem-unif-lip} the Lipschitz constants of $(\tilde{\gamma}_n)_n$ are bounded and since $\tilde{\gamma}_n(0)=p$ for all 
$n$, Theorem \ref{thm-cq-1.6} gives a $\ttc$-convergent subsequence, say $\tilde{\gamma}_{n_k}\to\tilde{\gamma}$.

Lemma \ref{lem-phi-cont} gives that $\gamma_{n_k} = \Phi(\tilde{\gamma}_{n_k})\to \Phi(\tilde{\gamma})$. The 
topology $\tau$ is Hausdorff, hence $\Phi(\tilde{\gamma})=\lambda$. Moreover, $d$ does not depend on the parametrization of 
the curves, so for $t\in[0,1]$ we get that $d^h(\tilde{\gamma}(t),\Gamma_{n_k})\leq 
d^h(\tilde{\gamma}(t),\tilde{\gamma}_{n_k}(t))\leq\rho(\tilde{\gamma},\tilde{\gamma}_{n_k})$ and analogously 
$d^h(\tilde{\gamma}_{n_k}(t),\Lambda)\leq\rho(\tilde{\gamma},\tilde{\gamma}_{n_k})$. Finally, this gives 
$d(\gamma_{n_k},\lambda)\to 0$ --- a contradiction to $d(\gamma_n,\lambda)\geq \ep$ for all $n\in\N, n\geq 1$.
\end{pr}

Actually, in all of the above (and results derived later on), where we used the non-total imprisonment condition, 
it would have sufficed to have the following condition on the space-time: For all $p,q\in M$, for every compact set $K$ with 
$p,q\in K$, there is a constant $C>0$ such that for all causal curves $\gamma$ from $p$ to $q$ contained in $K$ we have that 
$L^h(\gamma)\leq C$. This condition is weaker than non-totally imprisoning and stronger than causality, but at the moment it 
is unclear if it is strictly weaker or strictly stronger, respectively. Thus we opted for using an established causality 
condition, i.e., non-total imprisonment.

\begin{rem}\label{rem-top-sum}
In summary, when it comes to compactness (in a non-totally imprisoning spacetime), these results allow us to choose to work 
with unparametrized causal curves (i.e., on $\Cpq$) or causal curves with a fixed parametrization on $[0,1]$ (i.e., on 
$\Cpqt$).

Moreover, since $\tau$ and $\tau'$ agree for non-totally imprisoning spacetimes, the topology induced by the Hausdorff metric 
does not depend on the choice of the Riemannian background metric $h$ (the compact-open topology does not depend on the 
metric on the target space anyway, if the domain space is compact, see \cite[Theorem 4.2.17, p.\ 263]{Eng:89}). 
\end{rem}

\section{Global hyperbolicity}\label{sec-gh}
We want to extend the usual notion of global hyperbolicity in terms of the compactness of the causal diamonds 
to continuous spacetimes and investigate which equivalences still hold. We start by giving the definition and note that we 
do not require $(M,g)$ to be strongly causal nor stably causal (as in \cite{CG:12, FS:12}). However, just causality (as in
the smooth case, cf.\ \cite{BS:07}) is not enough, since causality and the compactness of the causal diamonds imply the
closedness of all $J^\pm(p)$, but not the closedness of the causal relation $\leq$. The reason for this is that one would
need (at least in the standard proofs) convex neighborhoods and the push-up lemma (cf.\ \cite[Lemma 1.22, p.\ 14]{CG:12}),
which one does not have in general (\cite[Examples 1.11-1.13, p.\ 9ff.]{CG:12}). Thus we use non-total
imprisonment as the causality requirement, which has been used in place of causality in the definition of global 
hyperbolicity for the smooth case (cf.\ \cite[Section 3, p.\ 831]{Min:09a}). Furthermore, as mentioned at the end of section 
\ref{sec-top}, it would actually suffice to have something slightly weaker than non-totally imprisoning, i.e., the condition 
above Remark \ref{rem-top-sum}.

\begin{defi}
 The spacetime $(M,g)$ is called \emph{globally hyperbolic} if
\begin{enumerate}
 \item $(M,g)$ is non-totally imprisoning and
 \item for all $p,q\in M$ the set $J(p,q):=J^+(p)\cap J^-(q)$ is compact (with respect to the manifold topology).
\end{enumerate}
\end{defi}

At this point we want to relate this notion of global hyperbolicity to other definitions of
global hyperbolicity. We will use extensively the results of section \ref{sec-top}. In particular, that in a non-totally
imprisoning spacetime one has $\tau=\tau'$ and $\Cpq$ is compact if and only if $\Cpqt$ is. Moreover the topology on $\Cpq$
and $\Cpqt$ will always be $\tau$ and $\ttc$, respectively.

\bigskip
The following theorem establishes that also in the continuous case global hyperbolicity is equivalent to $\Cpq$
being compact and $(M,g)$ being non-totally imprisoning. Moreover it indicates that these two notions are closely related. 
\begin{thm}\label{thm-gh-cpq}
Let $(M,g)$ be non-totally imprisoning. Then $(M,g)$ is globally hyperbolic if and only if for all $p,q\in
M$, $\Cpq$ is compact.
\end{thm}
\begin{pr}
Assume $(M,g)$ is globally hyperbolic. By Corollary \ref{cor-phi-proper} we can equivalently work with $\Cpqt$. So let
$(\gamma_n)_n$ be a sequence in $\Cpqt$ (showing sequential compactness suffices, since $\ttc$ is metrizable), then
the set $A:=\bigcup_{n\in\N}\Gamma_n$ is contained in the compact set $J(p,q)$ by definition. Therefore by Lemma
\ref{lem-unif-lip} the Lipschitz constants of the $\gamma_n$'s are uniformly bounded. Moreover since $\gamma_n(0)=p$ for all
$n\in\N$, there is a curve $\gamma\in\Cpqt$ such that a subsequence of $(\gamma_n)_n$ converges to $\gamma$ with respect to
$\ttc$ by Theorem \ref{thm-cq-1.6}. Therefore $\Cpqt$ is compact, and so $\Cpq$ is compact, too as the continuous image of
the compact set under the projection $\Phi$.

Let $p,q\in M$ and assume that $\Cpq$ is compact. Let $(x_n)_n$ be a sequence in $J(p,q)$. Consequently there is a sequence
$(\gamma_n)_n$ in $\Cpqt$ such that $x_n\in\Gamma_n$ for all $n$. Since $\Cpqt$ is compact (by Corollary
\ref{cor-phi-proper}) there is a subsequence $(\gamma_{n_k})_k$ of $(\gamma_n)_n$ and a causal curve $\gamma\in\Cpqt$ such
that $\gamma_{n_k}\to\gamma$ with respect to $\ttc$. Let $(t_k)_k$ be in $[0,1]$ such that $\gamma_{n_k}(t_k)=x_{n_k}$ for
all $k\in\N$. Without loss of generality (by passing to a subsequence) we can assume that $t_k\to t^*\in[0,1]$, then
$d^h(x_{n_k},\gamma(t^*)) \leq d^h(\gamma_{n_k}(t_k),\gamma(t_k))+d^h(\gamma(t_k),\gamma(t^*))$. Here both terms go to zero
since the first is bounded by $\rho(\gamma_{n_k},\gamma)$ and the second one because $t_k\to t^*$ and $\gamma$ is continuous.
Summing up, we conclude that $x_{n_k}=\gamma_{n_k}(t_k)\to\gamma(t^*)\in\Gamma\subseteq J(p,q)$. 
\end{pr}

Also in the continuous case, global hyperbolicity implies causal simplicity.
\begin{prop}\label{prop-gh-caus-sim}
 Let $(M,g)$ be globally hyperbolic, then $J^\pm$ is closed, i.e., if $p_n\leq q_n$ for all $n\in\N$ and $p_n\to p$, 
$q_n\to q$ then $p\leq q$. So $(M,g)$ is \emph{causally simple}.
\end{prop}
\begin{pr}
 Let $p_n\leq q_n$ for all $n\in\N$, $p_n\to p$ and $q_n\to q$. Let $p'\in \check{I}^-(p)$ and $q'\in \check{I}^+(q)$.
Then $p\in \check{I}^+(p')$ and $q\in\check{I}^-(q')$, which are open and hence, contain all but finitely many of the 
$p_n$'s and $q_n$'s, respectively. So for $n\geq n_0$ we have that $(q_n)_{n\geq n_0}$, $(p_n)_{n\geq n_0}$ are in the 
compact set $J^+(p')\cap J^-(q')$, and hence this yields a sequence $(\gamma_n)_{n\geq n_0}$ in 
$\tilde{C}(p',q')$ through $p_n$ and $q_n$. By assumption, $\tilde{C}(p',q')$ is compact, so there is a subsequence of
$(\gamma_n)_{n\geq n_0}$ that converges uniformly to a causal curve $\gamma$ from $p'$ to $q'$. Let $(t_k)_k$, $(s_k)_k$ be
in $[0,1]$ such that for all $k\in\N$: $\gamma_{n_k}(t_k)=p_{n_k}$ and $\gamma_{n_k}(s_k)=q_{n_k}$. Since $p_n\leq q_n$ for
all $n\in\N$, we have $t_k\leq s_k$ for all $k\in\N$. Without loss of generality we can assume that $t_k\to t^*\in[0,1]$ and
$s_k\to s^*\in[0,1]$, then clearly $\gamma(t^*)=p$ and $\gamma(s^*)= q$. It remains to show that $p\leq q$. Assume the
contrary, i.e., $t^*>s^*$, then there is an $\ep>0$ with $t^*-\ep > s^*$ and there is a $k_0\in\N$ such that for all $k\geq
k_0$: $t_k\in(t^*-\frac{\ep}{2}, t^*+\frac{\ep}{2})$ and $s_k\in(s^*-\frac{\ep}{2}, s^*+\frac{\ep}{2})$. However, for $k\geq
k_0$ this implies $s_k< s^*+\frac{\ep}{2} < t^*-\frac{\ep}{2} < t_k$ --- a contradiction to $t_k\leq s_k$ for all $k\in\N$.
\end{pr}
This, in particular, implies that $J^\pm(p)$ is closed for all $p\in M$.
\bigskip

An immediate consequence is the following result.
\begin{cor}\label{cor-comp}
 A non-totally imprisoning spacetime $(M,g)$ is globally hyperbolic if and only if $J^+(K_1)\cap J^-(K_2)$ is compact for all
compacta $K_1, K_2 \subseteq M$.
\end{cor}
\begin{pr}
 The reverse implication is trivial. So let $(M,g)$ be globally hyperbolic and $K_1, K_2\subseteq M$ be compact. Let 
$(p_n)_n$ be a sequence in $J^+(K_1)\cap J^-(K_2)$. Hence there are sequences $(r_n)_n$ in $K_1$, $(q_n)_n$ in $K_2$ such 
that for all $n\in\N$: $r_n\leq p_n \leq q_n$. By the compactness of $K_1$ and $K_2$ we can without loss of generality assume 
that $r_n\to r\in K_1$ and $q_n\to q\in K_2$ for $n\to\infty$. Now let $r'\in \check{I}^-(r)$ and $q'\in 
\check{I}^+(q)$. Then, since $\check{I}^+(r')$ and $\check{I}^-(q')$ are open and contain $r$ and $q$, 
respectively, we get that there is an $n_0\in\N$ such that for all $n\geq n_0$: $r_n\in \check{I}^+(r')$ and $q_n\in 
\check{I}^-(q')$. From $r_n\leq p_n\leq q_n$ for all $n\in\N$ it follows that for $n\geq n_0$: $p_n\in J(r',q')$, which 
is compact by assumption. Hence there is a subsequence of $(p_n)_n$, which converges and moreover its limit $p$ has to 
satisfy $r\leq p \leq q$ (by Proposition \ref{prop-gh-caus-sim}), hence $p\in J^+(K_1)\cap J^-(K_2)$. 
\end{pr}

A globally hyperbolic spacetime is non-partially imprisoning (cf.\ Definition \ref{def-imp}).
\begin{prop}\label{prop-gh-non-impr}
 Let $(M,g)$ be globally hyperbolic and $K\subseteq M$ compact. Let $\gamma$ be a future directed future inextendible causal
curve starting at $p:=\gamma(0)\in K$, then there is a $t_0>0$ such that for all $t\geq t_0: 
\gamma(t)\notin K$. An analogous statement holds for past inextendible curves.
\end{prop}
\begin{pr}
 Let $\gamma$ be defined on $[a,b)$ and assume the contrary, i.e., that there is a sequence $(t_k)_k$ in
$[a,b)$ converging to $b$ with $t_{2k}\in K$ and $t_{2k+1}\notin K$ for all $k\in \N$. This implies that
$\Gamma\subseteq J^+(K)\cap J^-(K)$, which is compact by Corollary \ref{cor-comp}, since for every $t\in[a,b)$ there is a
$t_{2k}>t$ and hence $\gamma(t)\in J^-(\gamma(t_{2k}))\subseteq J^-(K)$ --- a contradiction to $(M,g)$ being non-totally
imprisoning.
\end{pr}

In \cite{SW:96} Sorkin and Woolgar introduce the notion of $K$-causality, where $K$ is the smallest relation containing
$I^\pm$ that is transitive and closed, for a spacetime with a continuous metric. They use this concept to define a notion of
global hyperbolicity, called $K$-global hyperbolicity for convenience here, in terms of the compactness of the $K$-causal
diamonds and $K$-causality. Then they show that $K$-global hyperbolicity is equivalent to the compactness of causal curves
connecting two points with respect to the Vietoris topology. Moreover they show that for metrics in $\Con^2$ $K$-global
hyperbolicity agrees with the usual notion of global hyperbolicity. This seems quite close to our notion of global
hyperbolicity for continuous metrics, in particular, taking into account that the Vietoris topology (defined on all non-empty
closed subsets of a topological space) agrees with the Hausdorff topology on compact sets of a metric space (cf.\
\cite[Theorem 4.1, p.\ 80f.]{Fil:98}). However, Sorkin and Woolgar use piecewise $\Con^1$ curves for their definition of
causal curves, where we use locally Lipschitz continuous curves to have better limit curve theorems. So it remains open
whether these notions are equivalent for continuous metrics. It is known, however, that they agree for metrics in
$\Con^{0,1}$, i.e., metrics
that are locally Lipschitz continuous, by \cite[Proposition 1.21, p.\ 14]{CG:12} (in fact, they agree for \emph{causally
plain} spacetimes, cf.\ \cite[Definition 1.16, p.\ 13]{CG:12}). The question of compatibility comes down to $I^\pm =
I^\pm_{\mathrm{pcw. }\Con^1}$, which is still open in this regularity (i.e., for $g\in\Con$).

\section{Stability of global hyperbolicity}\label{sec-stab}
It is well-known in the smooth case that global hyperbolicity is stable, i.e., that for $(M,g)$ globally hyperbolic, 
there is a $g'\succ g$ such that $(M,g')$ is globally hyperbolic too (see \cite[Section 6, p.\ 447ff.]{Ger:70} or 
\cite[Theorem 2.6, p.\ 4f.]{BM:08}). Quite recently, it was shown that global hyperbolicity is also stable for so-called
\emph{continuous cone structures} in \cite{FS:12}, of which the forward lightcones (in a spacetime with a continuous
Lorentzian metric) are an example. However, the notion of global hyperbolicity in \cite{FS:12} requires stable causality, 
so we cannot use this result directly and, moreover, we then would have to show that there is a Lorentzian metric which has 
lightcones in between the original metric and this globally hyperbolic cone structure. So we will establish that global 
hyperbolicity implies stable causality and then give a self-contained proof of the stability of global hyperbolicity instead 
since the second step is the same anyways.
\bigskip

In the construction below we have to take the convex combination of two Lorentzian metrics which, in general, need
not be a \hyphenation{Lo-rentz-ian} metric, as the following simple
example shows. Set $g_1:= -dt^2 + dx^2$, $g_2:=dt^2 - dx^2$ on $\R^2$, then $\frac{1}{2}g_1 + \frac{1}{2}g_2 = 0$. If the
metrics are however such that one has lightcones smaller than the other, then any convex combination is again a Lorentzian
metric.
\begin{lem}\label{lem-conv-comb}
 Let $g_1\prec g_2$ be two continuous metrics on $M$, such that $(M,g_1)$ is time-orientable. Then any continuous
convex combination $g$ of $g_1$ and $g_2$ is again a continuous Lorentzian metric and $(M,g)$ is also time-orientable.
\end{lem}
\begin{pr}
 Let $X$ be a continuous $g_1$-timelike vector field (which exists by the time-orientability of $(M,g_1)$). Let $\chi\colon 
M\rightarrow [0,1]$ be continuous, then $g:= \chi\, g_1 + (1-\chi)\, g_2$ is again a continuous, symmetric $(0,2)$-tensor. 
Let $p\in M$ and define in a coordinate chart around $p$ a Lorentzian frame $e_1:=X, e_2,\ldots, e_n$ for $g_2$, where 
$e_2,\ldots,e_n$ are $g_2$-spacelike ($X$ is also $g_2$-timelike). Because $e_2,\ldots, e_n$ are also $g_1$-spacelike, this 
is also a Lorentzian frame for $g_1$. Furthermore $X$ is also $g$-timelike and $e_2,\ldots,e_n$ are $g$-spacelike, so this is 
a Lorentzian frame for $g$. This shows that $g$ is a Lorentzian metric.
\end{pr}

\begin{defi}
 A spacetime $(M,g)$ is called stably causal if there is a metric $\hat{g}\succ g$ such that $(M,\hat{g})$ is causal.
\end{defi}

The following theorem uses ideas of \cite[Theorem 3, p.\ 807ff.]{Min:09b} and \cite[Lemma 2, p.\ 241f.]{Min:09c}, where 
related results are proved for smooth metrics.
\begin{thm}\label{thm-gh-stab-caus}
 Let $(M,g)$ be globally hyperbolic, then $(M,g)$ is stably causal. In fact, there is a smooth metric $g'\succ g$ such
that $(M,g')$ is non-totally imprisoning.
\end{thm}
\begin{pr}
 We will first establish that $(M,g)$ is \emph{compactly stably non-totally imprisoning} (cf.\ \cite[Definition 5, p.\
806]{Min:09b} for compactly stably causal), i.e., that for all relatively compact open sets $B\subseteq M$ there is a metric
$g'$ with $g'\succ g$ on $B$, $g'=g$ on $M\backslash B$ and $(M,g')$ non-totally imprisoning.

Assume $(M,g)$ to be not compactly stably non-totally imprisoning, i.e., there is a relatively compact open set $B\subseteq
M$ and a sequence of metrics $(g_n)_n$ such that for all $n\in\N$ $g_n\succeq g_{n+1}\succeq g$, $g_n\succ g_{n+1}\succ g$ on
$B$, $g_n=g$ on $M\backslash B$, $(M,g_n)$ is totally imprisoning and $g_n \to g$ locally uniformly. Consequently for all
$n\in\N$ there is a compact set $K_n$ and an inextendible $g_n$-causal curve $\gamma_n$ contained in $K_n$. Without loss of
generality we can assume that each $\gamma_n$ is future inextendible and parametrized with respect to $h$-arclength. Set
$A:=\overline{(B)_1}=\{x\in M: d^h(x,B)\leq 1\}$, then $A$ is compact and every $\gamma_n$ must intersect $A$, since
otherwise: $\Gamma_n\subseteq M\backslash(B)_{\frac{1}{2}}\cap K_n$, which is a compact subset of $M\backslash B$. So,
$\gamma_n$ is $g$-causal and contained in a compact set --- a contradiction to $(M,g)$ non-totally imprisoning.

Moreover, we can assume that $\gamma_n(0)=:p_n$ is the first point of intersection of $\gamma_n$ with $A$ and consider
$\gamma_n$ to be defined on $[0,\infty)$. Furthermore, because $A$ is compact we can without loss of generality assume that
$p_n\to p\in A$. If there was a subsequence $(\gamma_{n_k})_k$ of $(\gamma_n)_n$ such that $\Gamma_{n_k}\subseteq A$ for all
$k\in\N$, then by Theorem \ref{thm-cq-1.6} there would be an inextendible $g$-causal limit curve contained in $A$, again a
contradiction to $(M,g)$ non-totally imprisoning. Thus we are allowed to assume $\Gamma_n\nsubseteq A$ for all $n\in\N$.

Again, by Theorem \ref{thm-cq-1.6} there is a $g$-causal inextendible limit curve $\gamma$ with $\gamma(0)=p$ and we
parametrize $\gamma$ with respect to $h$-arclength, hence $\gamma\colon[0,\infty)$. Then, since $(M,g)$ is non-partially
imprisoning (Proposition \ref{prop-gh-non-impr}) we have $0<t^0:=\max\{s\in(0,\infty):\gamma(s)\in A\}<\infty$ ($t^0>0$,
since by the above, all $\gamma_n$ have to enter $(B)_{\frac{1}{2}}$).

Let $q:=\gamma(t^0)\in\partial{A}$. Let $(\gamma_{n_k})_k$ be a subsequence of $(\gamma_n)_n$ such that
$\gamma_{n_k}\to\gamma$ locally uniformly. Since $\gamma_{n_k}(t^0+1)\to\gamma(t^0+1)\notin A$ we can without loss
of generality assume that $\gamma_{n_k}(t^0+1)\notin A$ for all $k\in\N$. For $k\in\N$ let $(s_k^0,s_k^1)$ be the
largest open connected interval containing $t^0+1$ such that $\gamma_{n_k}((s_k^0,s_k^1))\subseteq M\backslash A$ and let
$q^0_k:=\gamma_{n_k}(s^0_k)$, $q_k^1:=\gamma_{n_k}(s^1_k)$. If $s_k^1=\infty$, then $\gamma_{n_k}|_{[s_k^0,\infty)}$ would
be contained in the compact set $K_{n_k}\cap M\backslash (B)_{\frac{1}{2}}$, so $\gamma_{n_k}$ is $g$-causal there and hence
this gives a contradiction. Consequently $s_k^1<\infty$ for all $k\in\N$ and $q_k^0,q_k^1\in\partial A$. Since $\partial A$
is compact we are able to assume that $q^1_k\to q^1\in\partial A$. Also, $\gamma_{n_k}|_{[s_k^0,s_k^1]}$ is $g$-causal
because $\gamma_{n_k}([s_k^0,s_k^1])\subseteq M\backslash B$ and $g_{n_k}=g$ on $M\backslash B$. Furthermore
$s_k^0\in[0,t^0+1]$ because $t^0+1\in(s_k^0,s_k^1)$ and so we can assume that $s_k^0\to s^0\in[0,t^0+1]$. Now we claim that
$s^0\leq t^0$, indeed, if $s^0>t^0$, then $\gamma_{n_k}(s_k^0)=q_k^0\to \gamma(s^0)\in M\backslash A$, but
$q_k^0\in\partial A$, which is closed --- a contradiction.

Let $(r_k)_k$ be any sequence in $(0,t^0+1)$ with $s_k^0< r_k < t^0+1$ that converges to $t^0$,
then $w_k:=\gamma_{n_k}(r_k)\in M\backslash A$ for all $k\in\N$ and $w_k\to \gamma(t^0)=q$. This yields that the
$g$-causal curves $\gamma_{n_k}|_{[r_k,s_k^1]}$ have endpoints $w_k\leq q_k^1$ and hence $q\leq q^1$ by the closedness of
the causal relation (Proposition \ref{prop-gh-caus-sim}). However, since $q^1\in\partial A\cap \Gamma$ and $q$ is the last
point of $\gamma$ in $A$ it follows that $q^1\leq q$. In fact, we have that $q< q^1$, since $t^0+1< s_k^1$ for 
all $k\in\N$ (because $t^0+1\in(s_k^0, s_k^1)$) and so $\gamma_{n_k}(t^0+1)<\gamma_{n_k}(s_k^1)=q_k^1$ for all $k\in\N$. 
Then, again by the closedness of the causal relation $\gamma(t^0+1)\leq q^1$, which implies $q=\gamma(t^0)< 
\gamma(t^0+1)\leq q^1$. Finally, this gives the contradiction $q<q^1\leq q$.

Moreover, the spacetime $(M,g')$ constructed above is globally hyperbolic. Indeed, $(M,g')$ is non-totally imprisoning and
let $p,q\in M$ with $p\in J_{g'}^-(q)$. Let $(\gamma_n)_n$ be a sequence of $g'$-causal curves from $p$ to $q$. Then for
all $n\in\N$ we have 
\begin{equation*}
\Gamma_n\subseteq 
(J^+(p)\cap J^-(q))\cup (J^+(p)\cap J^-(\overline{B})) \cup (J^+(\overline{B})\cap J^-(q)) \cup (J^+(\overline{B})\cap 
J^-(\overline{B}))\,,
\end{equation*}
which is compact by Corollary \ref{cor-comp}. So, by Lemma \ref{lem-unif-lip} there is uniform bound 
on the Lipschitz constants of the $\gamma_n$'s, and hence by Theorem \ref{thm-cq-1.6} there is a converging subsequence.
This shows that $\tilde{C}_{g'}(p,q)$ is $\ttc$-compact and so $(M,g)$ is globally hyperbolic by Theorem \ref{thm-gh-cpq}
(via Corollary \ref{cor-phi-proper}). 

Clearly, from the above follows that for compacta $C, C'\subseteq M$ with $C\subseteq (C')^\circ$, there is
a metric $g'\succeq g$ with $g'\succ g$ on $C$, $g'=g$ on $M\backslash C'$ and $(M,g')$ globally hyperbolic.

Fix $w\in M$ and set $B_n:=\overline{B^{d^h}_n(w)}$, then $B_n$ is compact since $h$ is complete. Let $g_2\succeq g$ be a 
metric such that $g_2\succ g$ on $B_2$, $g_2=g$ on $M\backslash B_3$ and $(M,g_2)$ is globally hyperbolic. Define the 
compacta $C_3:=B_3\backslash B_2^\circ$, $C_3':=B_4\backslash B_1^\circ$ and let $g_3\succeq g_2$ be a metric with $g_3\succ 
g_2$ on $C_3$, $g_3=g_2$ on $M\backslash C_3'$ and $(M,g_3)$ globally hyperbolic. Inductively this yields a sequence 
$(g_n)_n$ of metrics satisfying $g_n\succeq g_{n-1}$, $g_n\succ g_{n-1}$ on $C_n=B_n\backslash B_{n-1}^\circ$, $g_n=g_{n-1}$ 
on $M\backslash C_n'$, where $C_n'=B_{n+1}\backslash B_{n-2}^\circ$ and $(M,g_n)$ globally hyperbolic. We claim that 
if $x\in B_n$, then $g_k(x)$ is independent of $k$, for $k\geq n+3$. Indeed, let $k\geq n+3$, then $x\in B_n\subseteq 
B_{n+1}^\circ\subseteq B_{k-2}^\circ$ and hence $g_k(x)=g_{k-1}(x)$ because $g_k=g_{k-1}$ on $M\backslash 
C_k'=B_{k-2}^\circ\cup M\backslash B_{k+1}$. Continuing in this way we find that $g_k(x)=g_{n+3}(x)$.
This allows us to define $g''(x):= g_{n+3}(x)$ for $x\in B_{n}$, where $n\in\N$ is minimal such that $x\in B_{n}$. 
By the above, this assignment is well-defined and $g''$ is a continuous Lorentzian metric satisfying $g''\succeq g_n$ for 
all $n\in\N$ and $g''\succ g$. It remains to show that $(M,g'')$ is non-totally-imprisoning. To this end suppose that $(M,
g'')$ is totally imprisoning, i.e., there is an inextendible $g''$-causal curve $\gamma$ contained in some compact set $K$.
Thus there is an $n\in\N$ such that $K \subseteq B_{n}$ and since $g''=g_{n+3}$ on $B_{n}$, $\gamma$ is
$g_{n+3}$-causal. This is a contradiction to $(M,g_{n+3})$ being non-totally imprisoning. Finally, let $g\prec\hat{g}\prec
g''$ be a smooth metric, then $(M,\hat{g})$ is non-total imprisoning too.
\end{pr}
Note that, given the stable causality of $(M,g)$, it is easy to show that there is a smooth metric $\hat{g}\succ g$ such
that $(M,\hat{g})$ is non-totally imprisoning by using the smooth causal ladder (cf.\ \cite{MS:08}). Nevertheless, in the
above proof we needed to show that the spacetime constructed ($(M,g')$) by widening the lightcones on $B$ is non-totally
imprisoning too. Thus it was the same amount of work to show the non-total imprisonment property directly instead of just
showing causality of $(M,g')$.
\bigskip

As mentioned in the introduction of this section, we opted for giving a self-contained proof of the stability of global 
hyperbolicity. Nevertheless we use the idea of proof of \cite[Theorem 1.2, p.\ 334f.]{FS:12}, in fact the following lemma is 
the analog of \cite[Corollary 3.5, p.\ 330 and Lemma 4.2, p.\ 332]{FS:12} in our setting.
\begin{lem}\label{lem-stab-gh}
 Let $(M,g)$ be globally hyperbolic, let $K\subseteq M$ be compact and let $(g_n)_n$ be a sequence of smooth metrics with 
$g\prec g_{n+1}\prec g_n$ for all $n\in\N$ that converges locally uniformly to $g$.
\begin{enumerate}
 \item\label{lem-stab-gh-len} Then there are $n_0\in\N$ and $L>0$ such that for all $n\geq n_0$ and for every $g_n$-causal 
curve $\gamma$ contained in $K$ we have $L^h(\gamma)\leq L$.
 \item For every $\ep>0$ there is an $n_1\in\N$ such that for all $n\geq n_1$ and for every $g_n$-causal curve 
$\gamma\colon[a,b]\rightarrow K$ there exists a $g$-causal curve $\lambda\colon[a,b]\rightarrow K$ with 
$\sup_{t\in[a,b]}d^h(\gamma(t),\lambda(t))\leq \ep$.
\end{enumerate}
\end{lem}
\begin{pr}
 \begin{enumerate}
  \item By Theorem \ref{thm-gh-stab-caus} there is a smooth metric $\hat{g}\succ g$ such that $(M,\hat{g})$ is non-totally
imprisoning. Now let $n_0\in\N$ be such that for all $n\geq n_0$ we have $g_n\prec \hat{g}$ on $K$ by Lemma
\ref{lem-small-metri}. Then the claim follows from Lemma \ref{lem-unif-lip} for $(M,\hat{g})$.
  \item Assume to the contrary that there is an $\ep>0$ such that for every $n\in\N$ there is a $g_n$-causal curve 
$\gamma_n\colon[a,b]\rightarrow K$ with $\sup_{t\in[a,b]}d^h(\gamma_n(t),\lambda(t))> \ep$ (for every $g$-causal curve 
$\lambda\colon[a,b]\rightarrow K$). By point \ref{lem-stab-gh-len} above we have $L^h(\gamma_n)\leq L$ for some $L>0$ (for 
$n\geq n_0$). Moreover, since $\Gamma_n \subseteq K$ (for all $n\in\N$) we can without loss of generality assume that 
$(\gamma_n)_n$ accumulates at some point in $K$, hence by Theorem \ref{thm-cq-1.6} there is a $g$-causal curve 
$\lambda\colon[a,b]\rightarrow K$ and a subsequence of $(\gamma_n)_n$ that converges uniformly to $\lambda$, contradicting 
$\sup_{t\in[a,b]}d^h(\gamma_n(t),\lambda(t))> \ep$ (for all $n\in\N$).
 \end{enumerate}
\end{pr}

With these preparations we can proceed to prove the stability of global hyperbolicity, where we use the idea of 
\cite[Theorem 1.2, p.\ 334f.]{FS:12}.
\begin{thm}\label{thm-gh-stab}
 Let $(M,g)$ be globally hyperbolic, then there is a smooth metric $g'\succ g$ such that $(M,g')$ is globally hyperbolic.
\end{thm}
\begin{pr}
 Fix $w\in M$ arbitrarily and set $M_1:=\overline{B_1^{d^h}(w)}$. We will inductively define a compact exhaustion $\bigcup_n 
M_n = M$. So assume $M_{n-1}$ has been constructed, then because $A_n:=\overline{(M_{n-1})_1} = \{x\in M: d^h(x,M_{n-1})\leq 
1\}$ is compact, $K_n := J^+(A_n)\cap J^-(A_n)$ is compact too (by Corollary \ref{cor-comp}). Consequently $M_n := 
\overline{(K_n)_1}$ is compact and $\bigcup_n M_n = M$.

At this point we define another covering of $M$ derived from $(M_n)_n$. Set $N_n:= M_n$ for $n=1,2$ and 
$N_n:=M_n\backslash M_{n-2}^\circ$ for $n\geq 3$. Let $(g_k)_k$ be a sequence of smooth metrics that converges locally 
uniformly to $g$ and such that $g\prec g_{k+1}\prec g_k$ for all $k\in\N$. Then for $n\in\N$ we get from Lemma 
\ref{lem-stab-gh} the existence of some $k_0(n)\in\N$ and $L_n>0$ such that for all $k\geq k_0$ and for every $g_k$-causal 
curve $\gamma$ contained in $N_n$, we have $L^h(\gamma)\leq L_n$. Moreover for every $g_k$-causal curve 
$\gamma\colon[a,b]\rightarrow N_n$ there is a $g$-causal curve $\lambda\colon[a,b]\rightarrow N_n$ with 
$\sup_{t\in[a,b]}d^h(\gamma(t),\lambda(t))< 1$.

This construction yields a subsequence $(g_{k_0(n)})_n$ of $(g_k)_k$, which for simplicity we denote by $(g_n')_n$. Note 
that this subsequence can be chosen to be monotonically decreasing, as $k_0(n)$ can be chosen to be greater than $k_0(i)$ 
for $1\leq i<n$ and from the construction in the proof of Lemma \ref{lem-stab-gh} we know that we can choose $g_n'$ to be 
non-totally imprisoning on $N_{n}$.

Now we want to construct a metric $g''\succ g$ such that $g''\prec g'_n$ on $N_n$ (for $n\in\N$). Observe that $(M_n)_n$ is 
an exhaustion of $M$ satisfying $M_n\subseteq M_{n+1}^\circ$ for all $n\in\N$ and $(N_n)_n$ satisfies $N_i\cap 
N_j=\emptyset$ for all $i,j\in\N$ with $|i-j|\geq 2$. We set $O_n:= M_n^\circ$ for $n=1,2$ and $O_n:=M_n^\circ\backslash 
M_{n-2}$ for $n\geq 3$. Then $(O_n)_n$ is on open cover of $M$ satisfying $O_i\cap O_j=\emptyset$ for $i,j\in\N$ with 
$|i-j|\geq 2$, since $O_i\subseteq N_i$. Let $(\chi_n)_n$ be a smooth partition of unity subordinated to $(O_n)_n$. Setting 
$g'':=\sum_{n=1}^\infty \chi_n g_{n+1}'$ yields a smooth symmetric $(0,2)$-tensor. By construction we have $g''|_{N_n} = 
\chi_{n-1} g_n'|_{N_n} + \chi_{n} g_{n+1}'|_{N_n} + \chi_{n+1} g_{n+2}'|_{N_n}$.\\
Consequently $g''$ is Lorentzian metric by Lemma \ref{lem-conv-comb}, $g\prec g''$ and $g''\prec g_n'$ on $N_n$ for every 
$n\in\N$.

It remains to show that $(M,g'')$ is globally hyperbolic. Let $\gamma\colon[\tilde{a},\tilde{b}]\rightarrow M$ be 
$g''$-causal and let $n\in\N$ be minimal such that $\gamma(\tilde{a})=:~p$, $\gamma(\tilde{b})=:q\in M_n$. We claim that 
$\Gamma\subseteq M_{n+1}$. Indeed, if we assume to the contrary that there is a $n'\in\N$ with $p,q\in M_{n'}$ and 
$n_0:=\min\{k\in\N: \Gamma\subseteq M_k\}\geq n'+2$, then a segment of $\gamma$, say $\gamma_0$, is contained in 
$N_{n_0}=M_{n_0}\backslash M_{n_0-2}^\circ$, which intersects $M_{n_0}\backslash M_{n_0-1}$. Let $\gamma_0$ be defined on 
$[a,b]$ with $\gamma_0(a),\gamma_0(b)\in \partial M_{n_0-2}$. Since $g''\prec g_{n_0}$ on $N_{n_0}$, there is a $g$-causal 
curve $\lambda\colon[a,b]\rightarrow N_{n_0}$ with $\sup_{t\in[a,b]}d^h(\gamma_0(t),\lambda(t))< 1$. This implies that 
$\lambda(a),\lambda(b)\in\overline{(M_{n_0-2})_1} = A_{n_0-1}$ and hence $\Lambda\subseteq K_{n_0-1}= J^+(A_{n_0-1})\cap 
J^-(A_{n_0-1})$. Consequently, $\Gamma_0\subseteq \overline{(K_{n_0-1})_1} = M_{n_0-1}$ --- a contradiction to $\Gamma_0\cap 
M_{n_0}\backslash M_{n_0-1} \neq \emptyset$. Thus, $\Gamma\subseteq M_{n+1}$ for all $g''$-causal curves $\gamma$ with 
endpoints $p,q\in M_n$. We prove inductively that we can bound the $h$-lengths of all $g''$-causal curves in $M_{n+1}$. The 
claim is true for $M_1$, where the bound is $L_1$ (since $M_1=N_1$). Assume we have a bound $L_n'$ for $M_n$, then since 
$g''\prec g_{n+1}'$ on $N_{n+1}$ we get $L^h(\gamma)\leq L_{n+1}$ if $\Gamma\subseteq N_{n+1}$. If not, i.e., $\Gamma\cap 
M_{n-1}^\circ\neq \emptyset$, let $\gamma$ be defined on $[a,b]$ and set $a':=\min \gamma^{-1}(M_{n-1})$, $b':=\max 
\gamma^{-1}(M_{n-1})$. Then $\gamma([a,a']),\gamma([b',b])\subseteq N_n$, hence their $h$-length is bounded by $L_n$.
Setting $\hat{\gamma}:=\gamma|_{[a',b']}$ yields a $g''$-causal curve $\hat{\gamma}\colon[a',b']\rightarrow M$ that has 
endpoints in $M_{n-1}$, so by the above $\hat{\Gamma}\subseteq M_n$, hence $L^h(\hat{\gamma})\leq 2 L_{n} + 
L_n'=:L_{n+1}'$.

Finally this means that $(M,g'')$ is non-totally imprisoning by Lemma \ref{lem-unif-lip} and also globally
hyperbolic, since if $(\gamma_k)_k$ is a sequence of $g''$-causal curves from $p\in M_n$ to $q\in M_n$, then
$L^h(\gamma_k)\leq L_{n+1}'$ for all $k\in\N$, which gives a uniform bound on their Lipschitz constants, when parametrized on
$[0,1]$ proportional to $h$-arclength (Lemma \ref{lem-cpq-equiv},\ref{lem-uniq-param}). Thus by Theorem \ref{thm-cq-1.6},
there is a subsequence of $(\gamma_k)_k$ that converges uniformly to a $g''$-causal curve connecting $p$ and $q$.
\end{pr}

\section{Cauchy hypersurfaces}\label{sec-chs}
In this section we establish that global hyperbolicity is equivalent to the existence of a Cauchy hypersurface. First we 
give the definition of a Cauchy hypersurface and then we explore some properties of Cauchy hypersurfaces. Note, that 
\cite{ONe:83} defines a Cauchy hypersurface in terms of timelike curves. However, in our setting it is easier to work with 
causal curves.

\begin{defi}
 A subset $S$ of $M$ is called a \emph{Cauchy hypersurface} if it is met exactly once by every inextendible causal curve.
\end{defi}

A Cauchy hypersurface is indeed a hypersurface.
\begin{prop}\label{prop-CHS-acausal}
 A Cauchy hypersurface is a closed acausal (i.e., two points on it cannot be connected by a causal curve) topological
hypersurface.
\end{prop}
\begin{pr}
 Let $S$ be a Cauchy hypersurface. Then clearly $S$ is acausal.

Let $\check{g}\prec g$ be a smooth metric, then obviously $S$ is a Cauchy hypersurface in $(M,\check{g})$, hence it is a
closed topological hypersurface by \cite[Lemma 14.29, p.\ 415f.]{ONe:83}.
\end{pr}

\begin{lem}\label{lem-CHS-causal}
 Let $(M,g)$ contain a Cauchy hypersurface $S$, then $(M,g)$ is causal.
\end{lem}
\begin{pr}
 Assume there is a closed causal curve $\gamma$ in $(M,g)$, then $\gamma$ is inextendible. Consequently either $\gamma$ meets
$S$ infinitely often or not at all. Both alternatives contradict the fact that $S$ is a Cauchy hypersurface.
\end{pr}

\begin{lem}\label{lem-CHS-disj}
Let $S$ be a Cauchy hypersurface. Then $M=I^-(S)\cup S \cup I^+(S)$, and these sets are disjoint.
\end{lem}
\begin{pr}
 Since $S$ is acausal and hence achronal, we have that $I^-(S)\cap I^+(S)=\emptyset$. Assume that
$q\in I^\pm(S)\cap S$, then there is a timelike curve from $q$ to $p\in S$. But by Lemma \ref{prop-CHS-acausal} $p$ has to be
equal to $q$ --- a contradiction to $(M,g)$ chronological.
\end{pr}
Since we also have $M=I^-_{\check{g}}(S)\cup S \cup I^+_{\check{g}}(S)$ (disjoint union) for any smooth metric 
$\check{g}\prec g$ (because $S$ is also a Cauchy hypersurface for $(M,\check{g})$) we get that
$M=\check{I}^-(S) \cup S \cup \check{I}^+(S)$ (again as a disjoint union). Furthermore,
$\check{I}^\pm(S) = I^\pm(S)$ and so $I^\pm(S)$ is open and $\partial I^\pm(S) = S$.
\bigskip

We need the following modification of \cite[Proposition 8.3.4, p.\ 202f.]{Wal:84}.
\begin{lem}\label{lem-wald-inext-intersects}
 Let $S$ be a Cauchy hypersurface and let $\lambda$ be an inextendible causal curve. Then $\lambda$ intersects
$\check{I}^-(S)$, $S$ and $\check{I}^+(S)$.
\end{lem}
\begin{pr}
 Assume that $\lambda$ does not intersect $\check{I}^+(S)$. Then without loss of generality we can assume that
$\lambda\colon (-\infty,\infty)\rightarrow M$ and $\lambda(0)=:p\in \check{I}^-(S)$. Moreover since $\lambda$ is inextendible
it has to meet $S$, say at $x:=\lambda(t_0)$. By assumption $\lambda((t_0,\infty))\subseteq \check{I}^-(S)$ and so
$q:=\lambda(t_1)\in \check{I}^-(S)$ for some $t_1>t_0$, hence there is a smooth metric $\check{g}\prec g$ and a future 
directed $\check{g}$-timelike curve $\gamma$ from $q$ to some point on $S$. Joining $\lambda|_{[0,t_1]}$ with $\gamma$, one 
obtains a future directed $g$-causal curve that meets $S$ twice --- a contradiction.
\end{pr}

As noted in section \ref{sec-top} strong causality implies non-partial imprisonment and hence non-total imprisonment. We
establish below that a spacetime containing a Cauchy hypersurface is strongly causal.
\begin{prop}\label{prop-CHS-str-causal}
 Let $(M,g)$ contain a Cauchy hypersurface $S$, then $(M,g)$ is strongly causal, i.e., for every $p\in M$ and every
neighborhood $U$ of $p$ there is a neighborhood $V$ of $p$ such that $V\subseteq U$ and every causal curve starting and
ending in $V$ is contained in $U$.
\end{prop}
\begin{pr}
 Assume that $(M,g)$ is not strongly causal, i.e., there is a $p\in M$, a neighborhood $U$ of $p$, a sequence of open
neighborhoods $(O_n)_n$ of $p$ such that $O_{n+1}\subseteq O_n\subseteq U$, $\bigcap_{n\in\N} O_n = \{p\}$ and there is a
sequence of causal curves $(\gamma_n)_n$ that start and end in $O_n$ but leave $U$. By the Limit curve theorem
(\cite[Theorem 3.1(2), p.\ 8ff.]{Min:08a} there is an inextendible limit curve $\gamma$ through p ($\gamma$ cannot be
closed by Lemma \ref{lem-CHS-causal}) and without loss of generality let $\gamma_n\to\gamma$ locally uniformly. Since
$M=\check{I}^-(S) \cup S \cup \check{I}^+(S)$, we have three cases. First, let $p\in\check{I}^+(S)$. Without loss of
generality we can assume that $O_n\subseteq \check{I}^+(S)$ for all $n\in\N$. Then each $\gamma_n$ cannot intersect
$\check{I}^-(S)$ since otherwise it would intersect $S$ twice. Consequently, $\gamma$ cannot enter $\check{I}^-(S)$
(otherwise infinitely many $\gamma_n$'s would). This gives a contradiction to Lemma \ref{lem-wald-inext-intersects}.
Analogously the second case, that is $p\in\check{I}^-(S)$, follows. The last case is $p\in S$. Each $\gamma_n$ has to
intersect $\check{I}^-(S)\cap M\backslash U$ or $\check{I}^+(S)\cap M\backslash U$. So without loss of generality (by passing
to a subsequence) we can assume that each $\gamma_n$ intersects $\check{I}^+(S)\cap M\backslash U$. Either $\gamma_n$ starts
in $O_n\cap\check{I}^-(S)$ or in $O_n\cap (S\cup \check{I}^+(S))$ and in both cases remains in $\check{I}^+(S)$. In the first
case let $p_n$ be the unique intersection point of $\gamma_n$ with $S$. By uniform convergence and closedness of $S$ we
conclude that $p_n\to p$ (if there a infinitely many such $\gamma_n$'s). In summary, we can consider $\gamma_n$ to start in
$S\cup\check{I}^+(S)$ and so, as above, there is an inextendible limit curve through $p$, which remains in the closed set
$\check{I}^+(S)\cup S$ --- once more a contradiction to Lemma \ref{lem-wald-inext-intersects}.
\end{pr}

At this point we are able to show that the existence of a Cauchy hypersurface implies global hyperbolicity. This is
now a direct generalization of \cite[Theorem 8.3.9, p.\ 206]{Wal:84}, since we prepared all the results we need for
continuous metrics. In particular, we use that since $\tau$ is metrizable, compactness is equivalent to sequential
compactness.
\begin{thm}\label{thm-chs-gh}
 Let $S$ be a Cauchy hypersurface in $(M,g)$, then $(M,g)$ is globally hyperbolic.
\end{thm}
\begin{pr}
As noted in Proposition \ref{prop-CHS-str-causal}, the spacetime $(M,g)$ is strongly causal, hence non-totally imprisoning.
By Theorem \ref{thm-gh-cpq} and Corollary \ref{cor-phi-proper} it suffices to show that $\Cpqt$ is compact with respect to
$\ttc$ for all $p,q\in M$. Let $p,q\in M$ and $(\lambda_n)_n$ be a sequence in $\Cpqt$.

We consider first the case that $p,q\in \check{I}^-(S)$. Removing the point $q$ from the spacetime yields a new
spacetime $\tilde{M}:=M\backslash\{q\}$, in which $(\lambda_n)_n$ is a sequence of future inextendible causal curves. Since
they start at $p$, by Theorem \ref{thm-cq-1.6} there is a future inextendible
causal curve $\lambda$ (in $\tilde{M})$ and a subsequence $(\lambda_{n_k})_k$ of $(\lambda_n)_n$ such that
$\lambda_{n_k}\to\lambda$ uniformly on compact sets. The curves $(\lambda_n)_n$ cannot enter $\check{I}^+(S)$ since otherwise
they would meet $S$ more than once. In $M$, however, $\lambda$ is either inextendible or it ends at $q$.
If $\lambda$ is inextendible, then because $\check{I}^+(S)$ is open and $(\lambda_{n_k})_k$ converges locally
uniformly to $\lambda$, the curve $\lambda$ cannot enter $\check{I}^+(S)$ (since otherwise the $\lambda_{n_k}$'s would),
which is a contradiction to Lemma \ref{lem-wald-inext-intersects}. In summary, $\lambda$ ends at $q$, so $\lambda\in\Cpqt$
and $\Cpqt$ is compact. Analogously one can show that $\Cpqt$ is compact, if $p,q\in \check{I}^+(S)$.

It remains to show compactness for the case $p\in \check{I}^-(S)$, $q\in \check{I}^+(S)$. As in the first
case there is a subsequence $(\lambda_{n_k})_k$ of $(\lambda_n)_n$ and a causal curve $\lambda$, which starts at $p$ and 
enters $\check{I}^+(S)$ such that $\lambda_{n_k}\to\lambda$ uniformly on compact sets. Reversing this method and viewing
$(\lambda_{n_k})_k$ as a sequence of past inextendible past directed causal curves starting at $q$ in
$\hat{M}:=M\backslash\{p\}$ yields a subsequence $(\hat{\lambda}_m)_m$ of $(\lambda_{n_k})_k$ and a past inextendible past
directed causal curve $\hat{\lambda}$ starting at $q$ such that $\hat{\lambda}_{m}\to\hat{\lambda}$ uniformly on compact sets
(as in the first case). Moreover $\hat{\lambda}$ has to enter $\check{I}^-(S)$ (again by Lemma
\ref{lem-wald-inext-intersects}). We claim that $\lambda\in\Cpqt$. Every $\lambda_{n_k}$
meets $S$ at some $t_k\in[0,1]$. By passing to yet another subsequence we can assume that $t_k\to t^*\in[0,1]$. If $\lambda$
would exist only on an interval $[0,s]$ with $s<t^*$, then $\lambda$ would be a limit of curves which do not intersect $S$
and thus $\lambda$ cannot leave the closed set $\check{I}^-(S)\cup S$ --- a contradiction since $\lambda$ has to enter
$\check{I}^+(S)$. Consequently, $\lambda$ has to be defined on $[0,t^*]$. By symmetry, $\hat{\lambda}$ is defined on
$[t^*,1]$ and so $\Lambda=\hat{\Lambda}$, which implies that $\lambda_{n_k} \to \lambda$ uniformly on $[0,1]$.
\end{pr}

This, in particular, implies that the interior of a Cauchy development is globally hyperbolic (as in the smooth case, cf.\ 
\cite[Theorem 14.38, p.\ 421f.]{ONe:83}). As noted in \cite[p.\ 23]{CG:12}, there are various difficulties when defining 
Cauchy developments based on timelike curves for arbitrary continuous metrics. Thus we follow \cite{CG:12} and define the 
Cauchy development of some set $S\subseteq M$ as
\begin{equation*}
 \D(S):=\D^+(S) \cup \D^-(S)\,,
\end{equation*}
where $\D^\pm(S):=\{p\in M:$ every past/future-directed past/future inextendible causal curve through $p$ meets $S\}$.

The following result does not assume $S$ to be a spacelike $\Con^1$-hypersurface as is assumed in \cite[Theorem 2.6, p.\ 
26f., Theorem 2.7, p.\ 28 and Theorem 2.9, p.\ 29]{CG:12}. Recall that a temporal function is a smooth function $f\colon 
M\rightarrow \R$ with past directed timelike gradient (\cite[p.\ 44]{BS:05}) and a Cauchy time function is a map $f\colon
M\rightarrow \R$, that is surjective, monotonically increasing along future directed causal curves and each $f^{-1}(\{t\})$
is a Cauchy hypersurface ($t\in\R$).
\begin{cor}
 Let $S\subseteq M$ be acausal, then the interior of $\D(S)$ is globally hyperbolic. Moreover, there is a smooth Cauchy 
time function on $\D(S)^\circ$.
\end{cor}
\begin{pr}
 Clearly, $S$ is a Cauchy hypersurface for the spacetime $(\D(S)^\circ, g|_{\D(S)^\circ})$, hence $\D(S)^\circ$ is globally 
hyperbolic by Theorem \ref{thm-chs-gh}. Furthermore, by Theorem \ref{thm-gh-stab} there is a smooth metric $g'\succ 
g|_{\D(S)^\circ}$ such that $(\D(S)^\circ, g')$ is globally hyperbolic. Thus by the smooth result (\cite[Theorem 
1.2, p.\ 44]{BS:05}) there is a Cauchy temporal function for $(\D(S)^\circ, g')$, which is also a smooth Cauchy time 
function for $(\D(S)^\circ, g|_{\D(S)^\circ})$.
\end{pr}

Conversely, global hyperbolicity implies the existence of a Cauchy hypersurface. This follows directly from the
stability of global hyperbolicity and the smooth theory.
\begin{thm}
 Let $(M,g)$ be globally hyperbolic, then there is a Cauchy hypersurface $S$ in $M$ and moreover, $M$ is homeomorphic to 
$\R\times S$.
\end{thm}

\begin{pr}
 By Theorem \ref{thm-gh-stab}, there is a smooth metric $g'\succ g$ such that $(M,g')$ is globally hyperbolic. Then by 
the smooth result (e.g.\ \cite[Theorem 11, p.\ 447]{Ger:70}) there is a Cauchy hypersurface $S$ for $(M,g')$, but $S$ is
also a Cauchy hypersurface for $(M,g)$ and $M$ is homeomorphic to $\R\times S$ (by \cite[Property 7, p.\ 444]{Ger:70} applied 
to $(M,g')$).
\end{pr}

Note that we can have, additionally to the above, that $(M,g')$ is isometric to $(\R\times S,-\beta dt^2 + h_t)$,
where $S$ is a smooth $g'$-spacelike Cauchy hypersurface, $\beta\colon\R\times S \rightarrow (0,\infty)$ smooth and $h_t$ 
is a family of $t$-dependent Riemannian metrics on the slices $\{t\}\times S$ (\cite[Theorem 1.1, p.\ 43f.]{BS:05}). 
Furthermore, there is also a Cauchy temporal function for $(M,g')$, which is also a smooth Cauchy time function for 
$(M,g)$ (\cite[Theorem 1.2, p.\ 44]{BS:05}).

\section{Maximal curves}\label{sec-max}
In the absence of an intrinsic (i.e., one that does not involve approximating smooth metrics) notion of geodesics in a 
spacetime with continuous metrics, we investigate \emph{maximal curves}. In sufficiently high regularity (i.e., 
$g\in\Con^{1,1}$) geodesics are locally maximizing. However in low regularity geodesics need not be maximizing, as can be 
seen from the classical example of Hartman and Wintner (\cite{HW:51}). They exhibit a two-dimensional Riemannian manifold 
with metric being $\Con^1$ and its derivatives being uniformly $\alpha$-H\"older continuous with $\alpha<1$ such that 
geodesics are not minimizing locally. So one cannot think of maximal curves as (pre-)geodesics in a regularity class below 
$\Con^{1,1}$.

\begin{defi}
 A curve $\lambda\in\Cpqt$ (or analogously for $\lambda\in\Cpq$) is called \emph{maximal} if there is a
neighborhood $U$ of $\lambda$ in $\Cpqt$ such that $L(\lambda)\geq L(\gamma)$ for all $\gamma\in U$. 
\end{defi}

We make use of the following lemma, whose proof is elementary.
\begin{lem}\label{lem-up-sc}
 Let $X$ be a topological space and $f,f_n\colon X\rightarrow [0,\infty)$ ($n\in\N$). Moreover, let $f_n$ be upper
semi-continuous, $f\leq f_n$ for all $n\in\N$ and $f_n\to f$ pointwise. Then $f$ is upper semi-continuous.
\end{lem}

At this point we are able to establish upper semi-continuity of the length functional with respect to the $\ttc$-topology on
$\Cpqt$. 
\begin{thm}\label{thm-L-up-sc}
 The length functional $L\colon\Cpqt\rightarrow [0,\infty)$ is upper semi-continuous (with respect to $\ttc$).
\end{thm}
\begin{pr}
 Let $(g_n)_n$ be a sequence of smooth metrics with $g_n\to g$ locally uniformly (given by Proposition \ref{prop-cg-2.1}) 
and $g\prec g_n$ for all $n\in\N$. Set $\delta_n:=\Delta(g_n,g)$. Then $\delta_n\searrow 0$ as $n\to\infty$. Note that 
$\Delta(g_n,g)=\delta_n$ implies that $-g(X,X)\leq -g_n(X,X)+\delta_n$ for all $X\in TM$ with $h(X,X)=1$.

Since the length functional for smooth metrics is upper semi-continuous (cf.\ \cite[Theorem 2.4, p.\ 54]{Min:08a}), $L_{g_n}$
is upper semi-continuous on $\tilde{C}_{g_n}(p,q)$, hence also on $\Cpqt\subseteq \tilde{C}_{g_n}(p,q)$ for all $n\in\N$.
Moreover, by dominated convergence we have that $L_{g_n}(\lambda)\to L(\lambda)$ for every $\lambda\in\Cpqt$. Let
$\lambda\in\Cpqt$, then 
\begin{equation*}
 L(\lambda)\leq \int_0^1{\sqrt{-g_n(\dot\lambda,\dot\lambda)+ C^2 \delta_n}} \leq L_{g_n}(\lambda)+C\sqrt{\delta_n}\,,
\end{equation*}
where $C=\Lip(\lambda)$ is a bound on $\sqrt{h(\dot\lambda,\dot\lambda)}$. Now we can apply Lemma
\ref{lem-up-sc}, since $L_{g_n}+C\sqrt{\delta_n}$ is upper
semi-continuous for all $n\in\N$ and converges to $L(\lambda)$ for $n\to\infty$.
\end{pr}

In the smooth case global hyperbolicity implies the existence of a maximizing geodesic joining any two causally related
points in the spacetime (cf.\ \cite[Proposition 14.19, p.\ 411]{ONe:83}). For continuous metrics we can show that an
analogous property holds for maximal curves.
\begin{prop}
 Let $\Cpqt$ be compact with respect to $\ttc$ (and non-empty), then there is a (globally) maximal causal curve connecting 
$p$ and $q$. In particular, if $(M,g)$ is globally hyperbolic, then there is a (globally) maximal curve connecting any two 
points in the spacetime which are causally related.
\end{prop}
\begin{pr}
 By Theorem \ref{thm-L-up-sc} the length functional $L$ is upper semi-continuous on $\Cpqt$. Since $\Cpqt$ is compact, $L$ 
attains a maximum there by \cite[Theorem 7.1, p.\ 29]{DiB:02}, this is then a globally maximal causal curve from $p$ to $q$.

If $(M,g)$ is globally hyperbolic, then $\Cpqt$ is compact for every pair of points $p,q$ in $M$ by Theorem \ref{thm-gh-cpq}
and Corollary \ref{cor-phi-proper}. Hence by the above $\Cpqt$ contains a globally maximal causal curve from $p$ to $q$, if
$q\in J^+(p)$ (i.e., $\Cpqt$ non-empty).
\end{pr}

An example of an extrinsic notion of geodesics are the so-called \emph{limit geodesics} of \cite[Definition 1.9(i), 
p.\ 7]{CG:12}. Now a natural question arises: Are limit geodesics maximizing? We can answer this question for a special 
class of limit geodesics.
\begin{prop}
 Let $\gamma\in\Cpqt$ be a limit geodesic such that:
\begin{enumerate}
 \item $g_n\to g$ locally uniformly with $g_n\succ g$ for all $n\in\N$ (given by \cite[Prop.\ 1.2, Eq.\ (1.8)]{CG:12}),
 \item for every $n\in\N$ there is a $g_n$-geodesic $\gamma_n$ such that $\gamma_n\to\gamma$ locally uniformly, and
 \item $\gamma_n$ is globally maximizing in $\tilde{C}_{g_n}(p,q)$ for all $n\in\N$ with respect to the corresponding
$\ttc$-topology.
\end{enumerate}
Then $\gamma$ is a (globally) maximal curve (on $\Cpqt$ with respect to $\ttc$).
\end{prop}

\begin{pr}
Note that by the proof of \cite[Prop.\ 1.2]{CG:12} we can without loss of generality assume that for all $n\in\N$, for all $X\in TM$ $g$-causal
\begin{equation}\label{eq-g-mon}
 -g(X,X)\leq -g_{n+1}(X,X)\leq -g_n(X,X) \,.
\end{equation}
Then for all $m\leq n$ we have $\gamma_n\in \tilde{C}_{g_n}(p,q)\subseteq\tilde{C}_{g_m}(p,q)$ and
\begin{equation}\label{eq-len-mon}
 L_{g_n}(\gamma_n)\leq L_{g_m}(\gamma_n)\,.
\end{equation}
Now let $\alpha\in\tilde{C}(p,q)$, $\eps>0$ and let $m\in\N$ such that
\begin{equation}\label{eq-len-app}
 |L_{g_m}(\gamma)-L(\gamma)|<\frac{\eps}{2}\,,
\end{equation}
which holds by dominated convergence. Moreover, as $L_{g_m}$ is upper semi-continuous, there is a neighborhood $U$ of $\gamma$ in $\tilde{C}_{g_m}(p,q)$ such that for all $\lambda\in U$ we have
\begin{equation}\label{eq-l-gm-usc}
 L_{g_m}(\lambda)\leq L_{g_m}(\gamma)+\frac{\eps}{2}\,.
\end{equation}
Let $n\geq m$ such that $\gamma_n\in U$ (as $\gamma_n\to \gamma$), then finally we get that
\begin{equation*}
 L(\alpha)\mathop{\leq}^{\eqref{eq-g-mon}}L_{g_n}(\alpha)\mathop{\leq}^{(*)}L_{g_n}(\gamma_n)\mathop{\leq}^{\eqref{eq-len-mon}}L_{g_m}(\gamma_n)\mathop{\leq}^{\eqref{eq-l-gm-usc}}L_{g_m}(\gamma)+\frac{\eps}{2}\mathop{\leq}^{\eqref{eq-len-app}} L(\gamma)+\eps\,,
\end{equation*}
where in $(*)$ we used that $\gamma_n$ is maximal in $\tilde{C}_{g_n}(p,q)$. This holds for all $\eps>0$ and so $\gamma$ is maximal in $\tilde{C}(p,q)$.
\end{pr}

\section{Conclusion}
We have shown that the definition of global hyperbolicity based on the compactness of causal diamonds and non-total
imprisonment can be extended to spacetimes with continuous metrics and that this notion is equivalent to compactness of
causal curves connecting two points and the existence of a Cauchy hypersurface. Moreover, it also implies causal simplicity
(closedness of the causal relation) and stable causality. Furthermore, it follows from this that in a globally hyperbolic
spacetime any two causally related points can be connected by a maximal curve.

It would be interesting to see if this notions differs from the concept of global hyperbolicity of Sorkin and Woolgar (based 
on $K$-causality, \cite{SW:96}). This question comes down to another open point, namely whether causality depends on the 
regularity of the curves used. In fact, the question is if $I^\pm$ in terms of locally Lipschitz continuous curves gives the 
same set as $I^\pm$ in terms of piecewise $\Con^1$ curves.

Additionally, the relation between global hyperbolicity and causal plainness (\cite[Definition 1.16, p.\ 13]{CG:12}) is 
open at the moment.

\section*{Acknowledgment}
The author is grateful to G\"unther H\"ormann, Michael Kunzinger and Roland Steinbauer for helpful discussions and detailed
feedback to drafts of this work and to Ettore Minguzzi for suggesting the proof of Lemma \ref{lem-unif-lip}.

This work was supported by project P25326 of the Austrian Science Fund.

\bibliographystyle{alphaabbr}
\bibliography{GlobalHyperbolicityForContinuousMetrics}
\addcontentsline{toc}{section}{References}

\end{document}